\newcommand{\xbD}{\Delta}

\newcommand{\xbO}{\Omega}
\newcommand{\xbP}{\Pi}

\newcommand{\xba}{\alpha}
\newcommand{\xbb}{\beta}

\newcommand{\xbe}{\in}
\newcommand{\xbf}{\phi}
\newcommand{\xbg}{\gamma}
\newcommand{\xbh}{\eta}

\newcommand{\xbl}{\lambda}
\newcommand{\xbm}{\mu}

\newcommand{\xbo}{\omega}

\newcommand{\xbq}{\psi}
\newcommand{\xbr}{\rho}
\newcommand{\xbs}{\sigma}

\newcommand{\xCK}{\times}

\newcommand{\xCN}{\neg}
\newcommand{\xCQ}{\emptyset}

\newcommand{\xCf}{\hspace{0.1em}}

\newcommand{\xcA}{\forall}
\newcommand{\xcB}{\stackrel{\subset}{\neq}}
\newcommand{\xcC}{\not\subseteq}

\newcommand{\xcE}{\exists}

\newcommand{\xcH}{\not\Rightarrow}
\newcommand{\xcI}{\not\Leftarrow}

\newcommand{\xcN}{\hspace{0.2em}\not\sim\hspace{-0.9em}\mid\hspace{0.8em}}

\newcommand{\xcS}{\bigcap}
\newcommand{\xcT}{\bot}

\newcommand{\xcV}{\bigcup}

\newcommand{\xcX}{\Box}

\newcommand{\xcb}{\subset}
\newcommand{\xcc}{\subseteq}
\newcommand{\xcd}{\supseteq}
\newcommand{\xce}{\not\in}
\newcommand{\xcf}{\supset}

\newcommand{\xch}{\Rightarrow}
\newcommand{\xci}{\Leftarrow}
\newcommand{\xcj}{\Leftrightarrow}
\newcommand{\xck}{\leq}
\newcommand{\xcl}{\vdash}
\newcommand{\xcm}{\models}
\newcommand{\xcn}{\hspace{0.2em}\sim\hspace{-0.9em}\mid\hspace{0.58em}}

\newcommand{\xco}{\vee}
\newcommand{\xcp}{\rightarrow}

\newcommand{\xcr}{\leftrightarrow}
\newcommand{\xcs}{\cap}
\newcommand{\xcu}{\wedge}
\newcommand{\xcv}{\cup}
\newcommand{\xcx}{\Diamond}

\newcommand{\xcz}{\Box}

\newcommand{\xDC}{\hspace{2em}}
\newcommand{\xDH}{\item }

\newcommand{\xdA}{\mbox{\boldmath$A$}}

\newcommand{\xdC}{\mbox{\boldmath$C$}}
\newcommand{\xdD}{\mbox{\boldmath$D$}}

\newcommand{\xdL}{\mbox{\boldmath$L$}}

\newcommand{\xdO}{\mbox{\boldmath$O$}}
\newcommand{\xdP}{\mbox{\boldmath$P$}}

\newcommand{\xdR}{\Re}

\newcommand{\xda}{{\cal A}}

\newcommand{\xdl}{{\cal L}}
\newcommand{\xdm}{{\cal M}}

\newcommand{\xdp}{{\cal P}}

\newcommand{\xdu}{{\cal U}}

\newcommand{\xdx}{{\cal X}}
\newcommand{\xdy}{{\cal Y}}
\newcommand{\xdz}{{\cal Z}}

\newcommand{\xEH}{ & }
\newcommand{\xEI}{\begin{itemize}}
\newcommand{\xEJ}{\end{itemize}}
\newcommand{\xEP}{ \\ }

\newcommand{\xEd}{\neq}
\newcommand{\xEh}{\begin{enumerate}}
\newcommand{\xEj}{\end{enumerate}}

\newcommand{\xeb}{\prec}

\newcommand{\xes}{\sqsubseteq}

\newcommand{\xex}{\lceil}

\newcommand{\xFO}{\parallel}

\newcommand{\xfA}{\mid}

\newcommand{\xfo}{\hookrightarrow}

\newcommand{\Xl}{\ldots}

\newcommand{\ol}{\overline}

\newcommand{\xssc}{\scriptsize}

\newcommand{\bl}{\begin{lemma} \rm}
\newcommand{\el}{\end{lemma}}
\newcommand{\br}{\begin{remark} \rm}
\newcommand{\er}{\end{remark}}
\newcommand{\be}{\begin{example} \rm}
\newcommand{\ee}{\end{example}}
\newcommand{\bco}{\begin{corollary} \rm}
\newcommand{\eco}{\end{corollary}}
\newcommand{\bc}{\begin{claim} \rm}
\newcommand{\ec}{\end{claim}}
\newcommand{\bfa}{\begin{fact} \rm}
\newcommand{\efa}{\end{fact}}
\newcommand{\bp}{\begin{proposition} \rm}
\newcommand{\ep}{\end{proposition}}
\newcommand{\bd}{\begin{definition} \rm}
\newcommand{\ed}{\end{definition}}
\newcommand{\bcs}{\begin{construction} \rm}
\newcommand{\ecs}{\end{construction}}
\newcommand{\bcd}{\begin{condition} \rm}
\newcommand{\ecd}{\end{condition}}
\newcommand{\bt}{\begin{theorem} \rm}
\newcommand{\et}{\end{theorem}}
\newcommand{\bn}{\begin{notation} \rm}
\newcommand{\en}{\end{notation}}
\newcommand{\bfi}{\begin{bild} \rm}
\newcommand{\efi}{\end{bild}}
\newcommand{\bsta}{\begin{statement} \rm}
\newcommand{\esta}{\end{statement}}
\newcommand{\bcom}{\begin{comment} \rm}
\newcommand{\ecom}{\end{comment}}
\newcommand{\bdia}{\begin{diagram} \rm}
\newcommand{\edia}{\end{diagram}}

\newcommand{\bfc}{\begin{figure}[htb] \begin{center}}
\newcommand{\efc}{\end{center} \end{figure}}

\documentclass{article}
\usepackage{amssymb,epic,eepic}

\title{
Reactive preferential structures and nonmonotonic consequence
}

\author{Dov M Gabbay
\thanks{
Dov.Gabbay@kcl.ac.uk, www.dcs.kcl.ac.uk/staff/dg
} \\
King's College, London
\thanks{
Department of Computer Science, King's College London, Strand,
London WC2R 2LS, UK
} \\ \\
Karl Schlechta
\thanks{
ks@cmi.univ-mrs.fr, karl.schlechta@web.de, http://www.cmi.univ-mrs.fr/ $\sim$ ks
} \\
Laboratoire d'Informatique Fondamentale de Marseille
\thanks{
UMR 6166, CNRS and Universit\'{e} de Provence,
Address: CMI, 39, rue Joliot-Curie, F-13453 Marseille Cedex 13, France
}
}

\begin{document}

\newtheorem{lemma}{Lemma}[section]
\newtheorem{theorem}[lemma]{Theorem}
\newtheorem{proposition}[lemma]{Proposition}
\newtheorem{corollary}[lemma]{Corollary}
\newtheorem{claim}[lemma]{Claim}
\newtheorem{fact}[lemma]{Fact}
\newtheorem{remark}[lemma]{Remark}
\newtheorem{definition}{Definition}[section]
\newtheorem{construction}{Construction}[section]
\newtheorem{condition}{Condition}[section]
\newtheorem{example}{Example}[section]
\newtheorem{notation}{Notation}[section]
\newtheorem{bild}{Figure}[section]
\newtheorem{comment}{Comment}[section]
\newtheorem{statement}{Statement}[section]
\newtheorem{diagram}{Diagram}[section]

\maketitle

\renewcommand{\labelenumi}
  {(\arabic{enumi})}
\renewcommand{\labelenumii}
  {(\arabic{enumi}.\arabic{enumii})}
\renewcommand{\labelenumiii}
  {(\arabic{enumi}.\arabic{enumii}.\arabic{enumiii})}
\renewcommand{\labelenumiv}
  {(\arabic{enumi}.\arabic{enumii}.\arabic{enumiii}.\arabic{enumiv})}

\setcounter{secnumdepth}{3}
\setcounter{tocdepth}{3}

\begin{abstract}

We introduce Information Bearing Relation Systems (IBRS) as an abstraction
of many logical systems. We then define a general semantics for IBRS, and
show that a special case of IBRS generalizes in a very natural way
preferential semantics and solves open representation problems for weak
logical systems. This is possible, as we can "break" the strong coherence
properties of preferential structures by higher arrows, i.e. arrows, which do
not go to points, but to arrows themselves.

\end{abstract}

\tableofcontents

%
%
%
\section{
Introduction
}
\subsection{
Overview
}

Our aim is twofold here. We want to

 \xEh

 \xDH introduce IBRS, an abstraction of many semantical structures
in nonclassical logics

 \xDH use IBRS as a generalization of preferential structures to give
weaker than preferential logics a semantics.

 \xEj

Many semantical structures for nonclassical logics use sets of models,
called possible worlds, together with a relation, e.g. accessibility in
the
case of Kripke models for modal logic, distance in the case of
Stalnaker-Lewis semantics for counterfactual conditionals, or comparison
of ``normality'' in the case of preferential structures. IBRS allow not only
relations between objects, i.e. possible worlds, but also ``higher''
relations, where e.g. an object can be in relation with a pair of objects
(which are in a ``lower'' relation), etc.

It is sometimes natural to see the (basic) relations as attacks:
if $a \xeb b,$ then a is considered more ``normal'' (or less $'' abnormal''
)$ than $b,$
or, we may say that a attacks b's normality, and we may write this as an
arrow
from a to $b.$ When we have now an attack from $c$ against this arrow,
i.e. the
relation $a \xeb b,$ it is natural to see now the original attack as
destroyed,
i.e. a does not attack $b$ any more. So $a \xeb b$ is still there, but
without
any effect. This is possible in an IBRS, see Diagram \ref{Diagram IBRS-Base}.

\vspace{30mm}

\begin{diagram}

\label{Diagram IBRS-Base}
\index{Diagram IBRS-Base}

\centering
\setlength{\unitlength}{1mm}
{\renewcommand{\dashlinestretch}{30}
\begin{picture}(110,150)(0,0)

\put(50,50){\circle{80}}
\put(50,50){\circle{40}}

\path(50,40)(50,60)
\path(49,57.3)(50,60)(51,57.3)
\put(50,39.2){\circle*{0.3}}
\put(50,60.8){\circle*{0.3}}
\put(50,38){\xssc{$a$}}
\put(50,62){\xssc{$b$}}

\path(20,50)(49,50)
\path(46.2,51)(49,50)(46.2,49)
\put(19.2,50){\circle*{0.3}}
\put(17,50){\xssc{$c$}}

\put(60,50){\xssc{$X$}}
\put(60,80){\xssc{$Y$}}

\put(30,2) {{\rm\bf Attacking an arrow}}

\end{picture}
}

\end{diagram}

\vspace{4mm}

The decisive property of preferential structures is the trivial fact that
if $a,b \xbe X \xcc Y,$ and $a \xeb b$ ``in $X'',$ then this will also
hold ''in $Y''.$ This creates
strong coherence properties. IBRS can break them: suppose $c \xbe Y-$X,
and
$c$ attacks $a \xeb b,$ then $a \xeb b$ is still there in $Y,$ but not
effective any more.
Thus, IBRS, or generalized preferential structures, as we will call this
variant of IBRS, can give a semantics to weaker than preferential logics.

Preferential structures themselves were introduced as abstractions
of Circumscription independently in  \cite{Sho87b} and  \cite{BS85}.
A precise definition of these structures is given below
in Definition \ref{Definition Pref-Str}.

In an abstract consideration of desirable properties a logic might have,
 \cite{Gab85}
examined rules a nonmonotonic consequence relation $ \xcn $ should
satisfy:

 \xEh

 \xDH $ \xbD, \xba \xcn \xba,$

 \xDH $ \xbD \xcn \xba $ $ \xch $ $( \xbD \xcn \xbb $ $ \xcj $ $ \xbD,
\xba \xcn \xbb ).$

 \xEj

Both, the semantic and the syntactic, approaches were connected
in  \cite{KLM90}, where a representation
theorem was proved, showing that the (stronger than $Gabbay' $s) system
$P$
corresponds to ``smooth'' preferential structures. System $P$ consists of

 \xEh

 \xDH $ \xbf \xcn \xbq,$ $ \xbf \xcn \xbq ' $ $ \xch $ $ \xbf \xcn \xbq
\xcu \xbq ',$
 \xDH $ \xbf \xcn \xbq,$ $ \xbf ' \xcn \xbq $ $ \xch $ $ \xbf \xco \xbf '
\xcn \xbq,$
 \xDH $ \xcl \xbf \xcr \xbf ' $ $ \xch $ $( \xbf \xcn \xbq $ $ \xcj $ $
\xbf ' \xcn \xbq ),$
 \xDH $ \xbf \xcn \xbq,$ $ \xcl \xbq \xcp \xbq ' $ $ \xch $ $ \xbf \xcn
\xbq ',$
 \xDH $ \xcl \xbf \xcp \xbf ' $ $ \xch $ $ \xbf \xcn \xbf ',$
 \xDH $ \xbf \xcn \xbq $ $ \xch $ $( \xbf \xcn \xbq ' $ $ \xcj $ $ \xbf
\xcu \xbq \xcn \xbq ' ).$

 \xEj

where $ \xcl $ is classical provability.

Details can be found in Definition \ref{Definition Log-Cond}.

To the authors' knowledge, a precise semantics for Gabbay's system was
still lacking.
We will give it here, applying the idea of IBRS, which allow to ``break''
the relations of preferential structures, and thus their strong
coherence conditions, as was already illustrated above.

We will also see that the usual definition of ``smoothness'' - for every
arrow
from a to $b,$ there has to be an arrow from $a' $ to $b,$ where $a' $ is
a minimal
element - is stronger
than needed, and a weaker variant, ``essential smoothness'' (see
Definition \ref{Definition Essentially-Smooth}) is sufficient
to generate the desired property of cumulativity on the logical side.
\subsection{
Introduction to IBRS
}
\index{Motivation IBRS}

The human agent in his daily activity has to deal with many situations
involving change. Chief among them are the following

 \xEh

 \xDH Common sense reasoning from available data. This involves
predication
of what unavailable data is supposed to be (nonmonotonic deduction)
but it is a defeasible prediction, geared towards immediate change.
This is formally known as nonmonotonic reasoning and is studied by
the nonmonotonic community.

 \xDH Belief revision, studied by a very large community. The agent
is unhappy with the totality of his beliefs which he finds internally
unacceptable (usually logically inconsistent but not necessarily so)
and needs to change/revise it.

 \xDH Receiving and updating his data, studied by the update community.

 \xDH Making morally correct decisions, studied by the deontic logic
community.

 \xDH Dealing with hypothetical and counterfactual situations. This
is studied by a large community of philosophers and AI researchers.

 \xDH Considering temporal future possibilities, this is covered by
modal and temporal logic.

 \xDH Dealing with properties that persist through time in the near
future and with reasoning that is constructive. This is covered by
intuitionistic logic.

 \xEj

All the above types of reasoning exist in the human mind and are used
continuously and coherently every hour of the day. The formal modelling
of these types is done by diverse communities which are largely distinct
with no significant communication or cooperation. The formal models
they use are very similar and arise from a more general theory, what
we might call:

``Reasoning with information bearing binary relations''.
\index{Definition IBRS}

\bd

$\hspace{0.01em}$


\label{Definition IBRS}

 \xEh
 \xDH An information bearing binary relation frame IBR, has the form
$(S, \xdR ),$ where $S$ is a non-empty set and $ \xdR $ is a subset of
$S,$ where $S$ is defined by induction as follows:

 \xEh

 \xDH $S_{0}=S$

 \xDH $S_{n+1}$ $=$ $S_{n} \xcv (S_{n} \xCK S_{n}).$

 \xDH $S$ $=$ $ \xcV \{S_{n}:n \xbe \xbo \}$

 \xEj

We call elements from $S$ points or nodes, and elements from $ \xdR $
arrows.
Given $(S, \xdR ),$ we also set $ \xdP ((S, \xdR )):=S,$ and $ \xdA ((S,
\xdR )):= \xdR.$

If $ \xba $ is an arrow, the origin and destination of $ \xba $ are
defined
as usual, and we write $ \xba:x \xcp y$ when $x$ is the origin, and $y$
the destination of the arrow $ \xba.$ We also write $o( \xba )$ and $d(
\xba )$ for
the origin and destination of $ \xba.$

 \xDH Let $Q$ be a set of atoms, and $ \xdL $ be a set of labels (usually
$\{0,1\}$ or $[0,1]).$ An information assignment $h$ on $(S, \xdR )$
is a function $h:Q \xCK \xdR \xcp \xdL.$

 \xDH An information bearing system IBRS, has the form
$(S, \xdR,h,Q, \xdL ),$ where $S,$ $ \xdR,$ $h,$ $Q,$ $ \xdL $ are as
above.

 \xEj

See Diagram \ref{Diagram IBRS} for an illustration.

\vspace{10mm}

\begin{diagram}

\centering
\setlength{\unitlength}{0.00083333in}
{\renewcommand{\dashlinestretch}{30}
\begin{picture}(4961,5004)(0,0)
\path(1511,1583)(611,3683)
\blacken\path(685.845,3584.520)(611.000,3683.000)(630.696,3560.885)(672.451,3539.613)(685.845,3584.520)
\path(1511,1583)(2411,3683)
\blacken\path(2391.304,3560.885)(2411.000,3683.000)(2336.155,3584.520)(2349.549,3539.613)(2391.304,3560.885)
\path(3311,1583)(4361,4133)
\blacken\path(4343.050,4010.616)(4361.000,4133.000)(4287.570,4033.461)(4301.603,3988.750)(4343.050,4010.616)
\path(3316,1574)(2416,3674)
\blacken\path(2490.845,3575.520)(2416.000,3674.000)(2435.696,3551.885)(2477.451,3530.613)(2490.845,3575.520)
\path(986,2783)(2621,2783)
\blacken\path(2501.000,2753.000)(2621.000,2783.000)(2501.000,2813.000)(2465.000,2783.000)(2501.000,2753.000)
\path(2486,2783)(2786,2783)
\blacken\path(2666.000,2753.000)(2786.000,2783.000)(2666.000,2813.000)(2630.000,2783.000)(2666.000,2753.000)
\path(3311,1583)(2051,2368)
\blacken\path(2168.714,2330.008)(2051.000,2368.000)(2136.987,2279.083)(2183.406,2285.509)(2168.714,2330.008)
\path(2166,2288)(1906,2458)
\blacken\path(2022.854,2417.439)(1906.000,2458.000)(1990.019,2367.221)(2036.567,2372.629)(2022.854,2417.439)

\put(1511,1358) {{\xssc $a$}}
\put(3311,1358) {{\xssc $d$}}
\put(3311,1058)  {{\xssc $(p,q)=(1,0)$}}
\put(1511,1058)  {{\xssc $(p,q)=(0,0)$}}
\put(2411,3758){{\xssc $c$}}
\put(4361,4433){{\xssc $(p,q)=(1,1)$}}
\put(4361,4208){{\xssc $e$}}
\put(2411,3983){{\xssc $(p,q)=(0,1)$}}
\put(611,3983) {{\xssc $(p,q)=(0,1)$}}
\put(611,3758) {{\xssc $b$}}
\put(1211,2883){{\xssc $(p,q)=(1,1)$}}
\put(260,2333) {{\xssc $(p,q)=(1,1)$}}
\put(2261,1583) {{\xssc $(p,q)=(1,1)$}}
\put(1286,3233){{\xssc $(p,q)=(1,1)$}}
\put(2711,3083){{\xssc $(p,q)=(1,1)$}}
\put(3836,2633){{\xssc $(p,q)=(1,1)$}}

\put(300,700)
{{\rm\bf
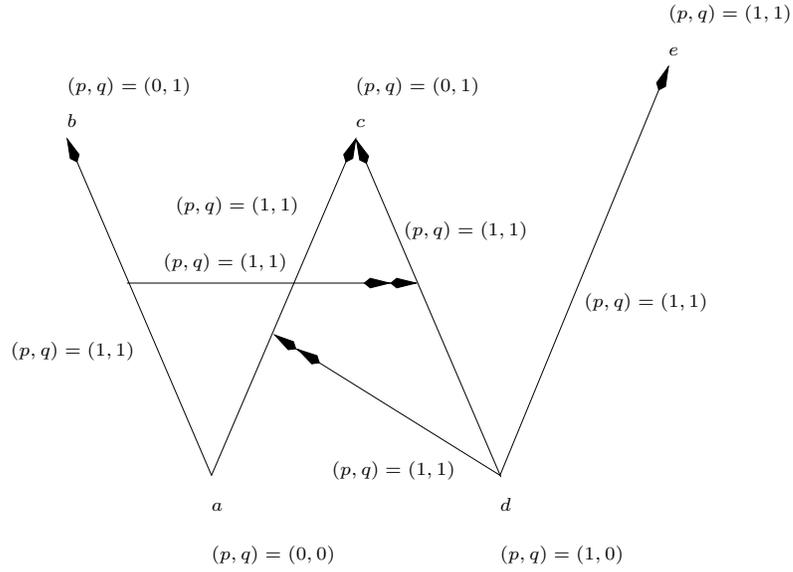
A simple example of an information bearing system.
}}

\end{picture}
}
\label{Diagram IBRS}
\index{Diagram IBRS}

\end{diagram}

We have here:
\[\begin{array}{l}
S =\{a,b,c,d,e\}.\\
\xdR = S \cup \{(a,b), (a,c), (d,c), (d,e)\} \cup \{((a,b), (d,c)),
(d,(a,c))\}.\\
Q = \{p,q\}
\end{array}
\]
The values of $h$ for $p$ and $q$ are as indicated in the figure. For
example $h(p,(d,(a,c))) =1$.

\vspace{4mm}

\index{Comment IBRS}

\ed

\bcom

$\hspace{0.01em}$


\label{Comment}

\label{Comment IBRS}

The elements in Figure Diagram \ref{Diagram IBRS} can be interpreted in
many ways,
depending on the area of application.

 \xEh

 \xDH The points in $S$ can be interpreted as possible worlds, or
as nodes in an argumentation network or nodes in a neural net or
states, etc.

 \xDH The direct arrows from nodes to nodes can be interpreted as
accessibility relation, attack or support arrows in an argumentation
networks, connection in a neural nets, a preferential ordering in
a nonmonotonic model, etc.

 \xDH The labels on the nodes and arrows can be interpreted as fuzzy
values in the accessibility relation or weights in the neural net
or strength of arguments and their attack in argumentation nets, or
distances in a counterfactual model, etc.

 \xDH The double arrows can be interpreted as feedback loops to nodes
or to connections.

 \xEj
\index{IBRS as abstraction}
\label{IBRS as abstraction}

\ecom

Thus, IBRS can be used as a source of information for various logics
based on the atoms in $Q.$ We now illustrate by listing several such
logics.
\subsubsection*{Modal Logic}

One can consider the figure as giving rise to two modal logic models.
One with actual world a and one with $d,$ these being the two minimal
points of the relation. Consider a language with $ \xcX q.$ how
do we evaluate $a \xcm \xcX q?$

The modal logic will have to give an algorithm for calculating the
values.

Say we choose algorithm $ \xda_{1}$ for $a \xcm \xcX q,$ namely:

[
$ \xda_{1}(a, \xcX q)=1$
]
iff for all $x \xbe S$ such that $a=x$ or $(a,x) \xbe \xdR $ we have
$h(q,x)=1.$

According to $ \xda_{1}$ we get that $ \xcX q$ is false at a. $ \xda_{1}$
gives rise to a $T-$modal logic. Note that the reflexivity is not
anchored at the relation $ \xdR $ of the network but in the algorithm
$ \xda_{1}$ in the way we evaluate. We say $(S, \xdR, \Xl.) \xcm \xcX $
$q$ iff $ \xcX q$ holds in all minimal points of $(S, \xdR ).$

For orderings without minimal points we may choose a subset of
distinguished
points.
\subsubsection*{Nonmonotonic Deduction}

We can ask whether $p \xcn q$ according to algorithm $ \xda_{2}$ defined
below. $ \xda_{2}$ says that $p \xcn q$ holds iff $q$ holds in all
minimal models of $p.$ Let us check the value of $ \xda_{2}$ in this
case:

Let $S_{p}=\{s \xbe S \xfA h(p,s)=1\}.$ Thus $S_{p}=\{d,e\}.$

The minimal points of $S_{p}$ are $\{d\}.$ Since $h(q,d)=0,$ we have that
$p \xcN q.$

Note that in the cases of modal logic and
nonmonotonic logic we ignored the arrows $(d,(a,c))$ (i.e. the double
arrow from $d$ to the arc $(a,c))$ and the $h$ values to arcs. These
values do not play a part in the traditional modal or nonmonotonic
logic. They do play a part in other logics. The attentive reader may
already suspect that we have her an opportunity for generalisation
of say nonmonotonic logic, by giving a role to arc annotations.
\subsubsection*{Argumentation Nets}

Here the nodes of $S$ are interpreted as arguments. The atoms $\{p,q\}$
can be interpreted as types of arguments and the arrows e.g. $(a,b) \xbe
\xdR $
as indicating that the argument a is attacking the argument $b.$

So, for example, let
\begin{quote}

a $=$ we must win votes.

$b$ $=$ death sentence for murderers.

$c$ $=$ We must allow abortion for teenagers

$d$ $=$ Bible forbids taking of life.

$q$ $=$ the argument is a social argument

$p$ $=$ the argument is a religious argument.

$(d,(a,c))$ $=$ there should be no connection between winning votes
and
abortion.

$((a,b),(d,c))$ $=$ If we attack the death sentence in order to win
votes then we must stress (attack) that there should be no connection
between religion (Bible) and social issues.
\end{quote}

Thus we have according to this model that supporting abortion can
lose votes. The argument for abortion is a social one and the argument
from the Bible against it is a religious one.

We can extract information from this IBRS using two algorithms. The
modal logic one can check whether for example every social argument
is attacked by a religious argument. The answer is no, since the social
argument $b$ is attacked only by a which is not a religious argument.

We can also use algorithm $ \xda_{3}$ (following Dung) to extract the
winning arguments of this system. The arguments a and $d$ are winning
since they are not attacked. $d$ attacks the connection between a
and $c$ (i.e. stops a attacking $c).$

The attack of a on $b$ is successful and so $b$ is out. However the
arc $(a,b)$ attacks the arc $(d,c).$ So $c$ is not attacked at all
as both arcs leading into it are successfully eliminated. So $c$
is in. $e$ is out because it is attacked by $d.$

So the winning arguments are $\{a,c,d\}$

In this model we ignore the annotations on arcs. To be consistent
in our mathematics we need to say that $h$ is a partial function on
$ \xdR.$ The best way is to give more specific definition on IBRS to
make it suitable for each logic.
\subsubsection*{Counterfactuals}

The traditional semantics for counterfactuals involves closeness of
worlds. The clauses $y \xcm p \xfo q,$ where $ \xfo $
is a counterfactual implication is that $q$ holds in all worlds $y' $
``near enough'' to $y$ in which $p$ holds. So if we interpret the
annotation on arcs as distances then we can define ``near'' as distance
$ \xck $ 2, we get: $a \xcm p \xfo q$ iff in all worlds
of $p-$distance $ \xck 2$ if $p$ holds so does $q.$ Note that the distance
depends on $p.$

In this case we get that $a \xcm p \xfo q$ holds. The
distance function can also use the arrows from arcs to arcs, etc.
There are many opportunities for generalisation in our IBRS set up.
\subsubsection*{Intuitionistic Persistence}

We can get an intuitionistic Kripke model out of this IBRS by letting,
for $t,s \xbe S,t \xbr_{0}s$ iff $t=s$ or $[tRs \xcu \xcA q \xbe Q(h(q,t)
\xck h(q,s))].$
We get that

[
$r_{0}=\{(y,y) \xfA y \xbe S\} \xcv \{(a,b),(a,c),(d,e)\}.$
]

Let $ \xbr $ be the transitive closure of $ \xbr_{0}.$ Algorithm $
\xda_{4}$
evaluates $p \xch q$ in this model, where $ \xch $ is intuitionistic
implication.

$ \xda_{4}:$ $p \xch q$ holds at the IBRS iff $p \xch q$ holds
intuitionistically
at every $ \xbr -$minimal point $of(S, \xbr ).$
\subsection{Purpose of this paper}

In this paper, we will not cover all these applications of the above
abstract
definition of IBRS, but will

(1) give an abstract and also a very concrete semantics to IBRS

(2) show that a special case of IBRS generalizes in a very natural way
preferential semantics and solves open representation problems for weak
logical systems. This is possible, as we can ``break'' the strong coherence
properties of preferential structures by higher arrows, i.e. arrows, which
do
not go to points, but to arrows themselves.

\section{
A semantics for IBRS
}

\label{Section Reac-Sem}
\index{Section Reac-Sem}

\subsection{
Introduction
}


\label{Section Reac-Sem-Intro}
\index{Section Reac-Sem-Intro}

(1) Nodes and arrows

As we may have counterarguments not only against nodes, but also against
arrows, they must be treated basically the same way, i.e. in some way
there has
to be a positive, but also a negative influence on both. So arrows cannot
just
be concatenation between the contents of nodes, or so.

We will differentiate between nodes and arrows by labelling arrows in
addition
with a time delay. We see nodes as situations, where the output is
computed
instantenously from the input, whereas arrows describe some ``force'' or
``mechanism'' which may need some time to ``compute'' the result from the
input.

Consequently, if $ \xba $ is an arrow, and $ \xbb $ an arrow pointing to $
\xba,$ then it
should point to the input of $ \xba,$ i.e. before the time lapse.
Conversely,
any arrow originating in $ \xba $ should originate after the time lapse.

Apart this distinction, we will treat nodes and arrows the same way, so
the
following discussion will apply to both - which we call just ``objects''.

(2) Defeasibility

The general idea is to code each object, say $X,$ by $I(X):U(X) \xcp
C(X):$ If $I(X)$
holds then, unless $U(X)$ holds, consequence $C(X)$ will hold. (We adopted
Reiter's
notation for defaults, as IBRS have common points with the former.)

The situation is slightly more complicated, as there can be several
counterarguments, so $U(X)$ really is an ``or''. Likewise, there can be
several
supporting arguments, so $I(X)$ also is an ``or''.

A counterargument must not always be an argument against a specific
supporting
argument, but it can be. Thus, we should admit both possibilties. As we
can use
arrows to arrows, the second case is easy to treat (as is the dual, a
supporting
argument can be against a specific counterargument). How do we treat the
case
of unspecific pro- and counterarguments? Probably the easiest way is to
adopt
Dung's idea: an object is in, if it has at least one support, and no
counterargument - see  \cite{Dun95}.
Of course, other possibilities may be adopted, counting, use
of labels, etc., but we just consider the simple case here.

(3) Labels

In the general case, objects stand for some kind of defeasible
transmission.
We may in some cases see labels as restricting this transmission to
certain
values. For instance, if the label is $p=1$ and $q=0,$ then the $p-$part
may be
transmitted and the $q-$part not.

Thus, a transmission with a label can sometimes be considered as a family
of
transmissions, which ones are active is indicated by the label.

\be

$\hspace{0.01em}$


\label{Example 2.1}

In fuzzy Kripke models, labels are elements of $[0,1].$ $p=0.5$ as label
for a node
$m' $ which stands for a fuzzy model means that the value of $p$ is 0.5.
$p=0.5$ as
label for an arrow from $m$ to $m' $ means that $p$ is transmitted with
value 0.5.
Thus, when we look from $m$ to $m',$ we see $p$ with value
$0.5*0.5=0.25.$ So, we have
$ \xcx p$ with value 0.25 at $m$ - if, e.g., $m,m' $ are the only models.

\ee

(4) Putting things together

If an arrow leaves an object, the object's output will be connected to the
(only) positive input of the arrow. (An arrow has no negative inputs from
objects it leaves.) If a positive arrow enters an object, it is connected
to
one of the positive inputs of the object, analogously for negative arrows
and
inputs.

When labels are present, they are transmitted through some operation.

\vspace{7mm}


\vspace{7mm}


\subsection{
Formal definition
}


\label{Section Reac-Sem-Def}
\index{Section Reac-Sem-Def}

\bd

$\hspace{0.01em}$


\label{Definition 2.1}

In the most general case, objects of IBRS have the form:
$(<I_{1},L_{1}>, \Xl,<I_{n},L_{n}>):(<U_{1},L'_{1}>, \Xl
,<U_{n},L'_{n}>),$ where the $L_{i},L'_{i}$ are labels
and the $I_{i},U_{i}$ might be just truth values, but can also be more
complicated,
a (possibly infinite) sequence of some values. Connected objects have,
of course, to have corresponding such sequences. In addition, the object
$X$ has a criterion for each input, whether it is valid or not (in the
simple
case,
this will just be the truth value $'' true'' ).$ If there is at least one
positive
valid input $I_{i},$ and no valid negative input $U_{i},$ then the output
$C(X)$ and its
label are calculated on the basis of the valid inputs and their labels. If
the
object is an arrow, this will take some time, $t,$ otherwise, this is
instantaneous.

\ed

Evaluating a diagram

An evaluation is relative to a fixed input, i.e. some objects will be
given
certain values, and the diagram is left to calculate the others. It may
well
be that it oscillates, i.e. shows a cyclic behaviour. This may be true for
a
subset of the diagram, or the whole diagram. If it is restricted to an
unimportant part, we might neglect this. Whether it oscillates or not can
also depend on the time delays of the arrows (see Example \ref{Example 2.2}).

We therefore define for a diagram $ \xbD $

$ \xba \xcn_{ \xbD } \xbb $ iff

(a) $ \xba $ is a (perhaps partial) input - where the other values are set
``not valid''

(b) $ \xbb $ is a (perhaps partial) output

(c) after some time, $ \xbb $ is stable, i.e. all still possible
oscillations
do not affect $ \xbb $

(d) the other possible input values do not matter, i.e. whatever the
input,
the result is the same.

In the cases examined here more closely, all input values will be defined.

\vspace{7mm}


\vspace{7mm}


\subsection{
A circuit semantics for simple IBRS without labels
}


\label{Section Reac-Sem-Circuit}
\index{Section Reac-Sem-Circuit}

It is standard to implement the usual logical connectives by electronic
circuits. These components are called gates. Circuits with feedback
sometimes
show undesirable behaviour when the initial conditions are not specified.
(When
we switch a circuit on, the outputs of the individual gates can have
arbitrary
values.) The technical
realization of these initial values shows the way to treat defaults. The
initial
values are set via resistors (in the order of 1 $k \xbO )$ between the
point in the
circuit we want to intialize and the desired tension (say 0 Volt for
false,
5 Volt for true). They are called pull-down or pull-up resistors (for
default
0 or 5 Volt). When a ``real'' result comes in, it will override the tension
applied via the resistor.

Closer inspection reveals that we have here a 3 level default situation:
The initial value will be the weakest, which can be overridden by any
``real''
signal, but a positive argument can be overridden by a negative one. Thus,
the biggest resistor will be for the initialization, the smaller one for
the
supporting arguments, and the negative arguments have full power.

Technical details will be left to the experts.

We give now an example which shows that the delays of the arrows can
matter.
In one situation, a stable state is reached, in another, the circuit
begins to
oscillate.

\be

$\hspace{0.01em}$


\label{Example 2.2}

(In engineering terms, this is a variant of a JK flip-flop with $R*S=0,$ a
circuit
with feedback.)

We have 8 measuring points.

$In1,In2$ are the overall input, $Out1,Out2$ the overall output,
$A1,A2,A3,A4$ are auxiliary internal points. All points can be true or
false.

The logical structure is as follows:

A1 $=$ $In1 \xcu Out1,$ A2 $=$ $In2 \xcu Out2,$

A3 $=$ $A1 \xco Out2,$ A4 $=$ $A2 \xco Out1,$

Out1 $=$ $ \xCN A3,$ Out2 $=$ $ \xCN A4.$

Thus, the circuit is symmetrical, with In1 corresponding to In2, A1 to A2,
A3 to A4, Out1 to Out2.

The input is held constant. See Diagram \ref{Diagram Gate-Sem}.

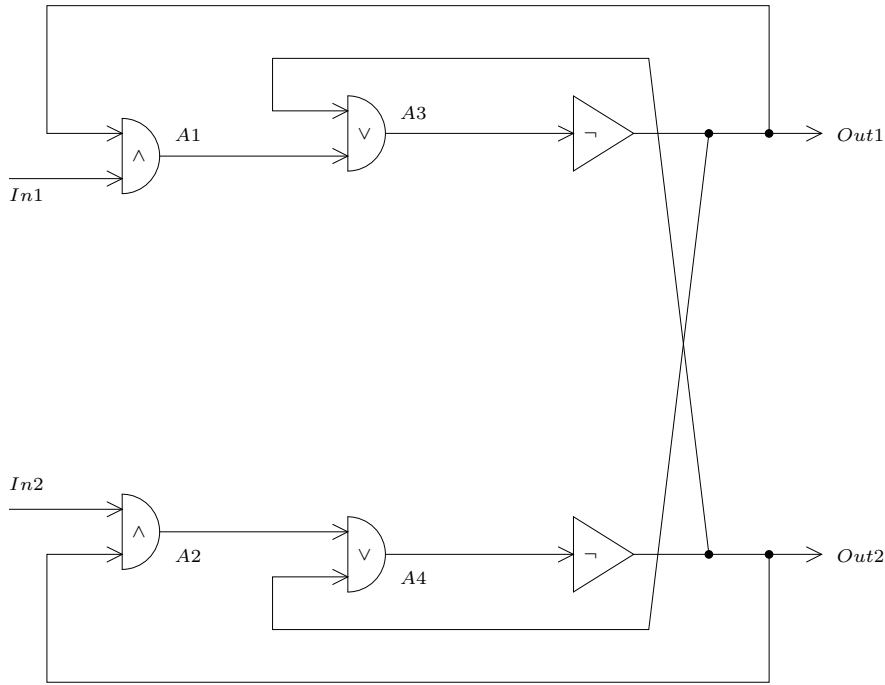
\vspace{10mm}

\begin{diagram}

\label{Diagram Gate-Sem}
\index{Diagram Gate-Sem}

\centering
\setlength{\unitlength}{1mm}
{\renewcommand{\dashlinestretch}{30}
\begin{picture}(150,170)(0,0)


\put(15,130){\arc{10}{-1.57}{1.57}}
\path(15,125)(15,135)
\put(16.3,129.3){\xssc{$\xcu$}}
\put(22,132){\xssc{$A1$}}
\path(13,132)(15,133)(13,134)
\path(13,126)(15,127)(13,128)

\put(45,133){\arc{10}{-1.57}{1.57}}
\path(45,128)(45,138)
\put(46.3,132.3){\xssc{$\xco$}}
\put(52,135){\xssc{$A3$}}
\path(43,135)(45,136)(43,137)
\path(43,129)(45,130)(43,131)

\path(75,128)(75,138)(83,133)(75,128)
\put(76.3,132.3){\xssc{$\xCN$}}
\path(73,132)(75,133)(73,134)

\path(0,127)(15,127)
\path(20,130)(45,130)
\path(50,133)(75,133)
\path(83,133)(108,133)
\path(106,132)(108,133)(106,134)
\put(93,133){\circle*{1}}
\put(101,133){\circle*{1}}

\put(0,124){\xssc{$In1$}}
\put(110,132){\xssc{$Out1$}}

\path(15,133)(5,133)(5,150)(101,150)(101,133)
\path(45,136)(35,136)(35,143)(85,143)(93,77)


\put(15,80){\arc{10}{-1.57}{1.57}}
\path(15,75)(15,85)
\put(16.3,79.3){\xssc{$\xcu$}}
\put(22,76){\xssc{$A2$}}
\path(13,82)(15,83)(13,84)
\path(13,76)(15,77)(13,78)

\put(45,77){\arc{10}{-1.57}{1.57}}
\path(45,72)(45,82)
\put(46.3,76.3){\xssc{$\xco$}}
\put(52,73){\xssc{$A4$}}
\path(43,79)(45,80)(43,81)
\path(43,73)(45,74)(43,75)

\path(75,72)(75,82)(83,77)(75,72)
\put(76.3,76.3){\xssc{$\xCN$}}
\path(73,78)(75,77)(73,76)

\path(0,83)(15,83)
\path(20,80)(45,80)
\path(50,77)(75,77)
\path(83,77)(108,77)
\path(106,76)(108,77)(106,78)
\put(93,77){\circle*{1}}
\put(101,77){\circle*{1}}

\put(0,85.5){\xssc{$In2$}}
\put(110,76){\xssc{$Out2$}}

\path(15,77)(5,77)(5,60)(101,60)(101,77)
\path(45,74)(35,74)(35,67)(85,67)(93,133)

\put(30,20) {{\rm\bf Gate Semantics}}

\end{picture}
}

\end{diagram}

\vspace{4mm}

\ee

We suppose that the output of the individual gates is present $n$ time
slices
after the input was present. $n$ will in the first circuit be equal to 1
for all
gates, in the second circuit equal to 1 for all but the AND gates, which
will
take 2 time slices. Thus, in both cases, e.g. Out1 at time $t$ will be the
negation of A3 at time $t-1.$ In the first case, A1 at time $t$ will be
the
conjunction of In1 and Out1 at time $t-1,$ and in the second case the
conjunction
of In1 and Out1 at time $t-2.$

We initialize In1 as true, all others as false. (The initial value of A3
and A4
does not matter, the behaviour is essentially the same for all such
values.)

The first circuit will oscillate with a period of 4, the second circuit
will go
to a stable state.

We have the following transition tables (time slice shown at left):

Circuit 1, $delay=1$ everywhere:

\begin{tabular}{lccccccccl}
   & In1 & In2 & A1 & A2 & A3 & A4 & Out1 & Out2 & \\
1: &   T &   F &  F &  F &  F &  F &    F &    F & \\
2: &   T &   F &  F &  F &  F &  F &    T &    T & \\
3: &   T &   F &  T &  F &  T &  T &    T &    T & \\
4: &   T &   F &  T &  F &  T &  T &    F &    F & \\
5: &   T &   F &  F &  F &  T &  F &    F &    F & oscillation starts \\
6: &   T &   F &  F &  F &  F &  F &    F &    T & \\
7: &   T &   F &  F &  F &  T &  F &    T &    T & \\
8: &   T &   F &  T &  F &  T &  T &    F &    T & \\
9: &   T &   F &  F &  F &  T &  F &    F &    F &
back to start of oscillation \\
\end{tabular}

Circuit 2, $delay=1$ everywhere, except for AND with $delay=2:$

(Thus, A1 and A2 are held at their intial value up to time 2, then they
are
calculated using the values of time $t-2.)$

\begin{tabular}{lccccccccl}
   & In1 & In2 & A1 & A2 & A3 & A4 & Out1 & Out2 & \\
1: &   T &   F &  F &  F &  F &  F &    F &    F & \\
2: &   T &   F &  F &  F &  F &  F &    T &    T & \\
3: &   T &   F &  F &  F &  T &  T &    T &    T & \\
4: &   T &   F &  T &  F &  T &  T &    F &    F & \\
5: &   T &   F &  T &  F &  T &  F &    F &    F & \\
6: &   T &   F &  F &  F &  T &  F &    F &    T & stable state reached \\
7: &   T &   F &  F &  F &  T &  F &    F &    T & \\
\end{tabular}

Note that state 6 of circuit 2 is also stable in circuit 1, but it is
never
reached in that circuit.

\vspace{7mm}


\vspace{7mm}


\section{
IBRS as generalized preferential structures
}


\subsection{
Introduction
}
\label{Section Log-Introduction}
\index{Definition Alg-Base}

\bd

$\hspace{0.01em}$


\label{Definition Alg-Base}

We use $ \xdp $ to denote the power set operator,
$ \xbP \{X_{i}:i \xbe I\}$ $:=$ $\{g:$ $g:I \xcp \xcV \{X_{i}:i \xbe I\},$
$ \xcA i \xbe I.g(i) \xbe X_{i}\}$ is the general cartesian
product, $card(X)$ shall denote the cardinality of $X,$ and $V$ the
set-theoretic
universe we work in - the class of all sets. Given a set of pairs $ \xdx
,$ and a
set $X,$ we denote by $ \xdx \xex X:=\{<x,i> \xbe \xdx:x \xbe X\}.$ When
the context is clear, we
will sometime simply write $X$ for $ \xdx \xex X.$

$A \xcc B$ will denote that $ \xCf A$ is a subset of $B$ or equal to $B,$
and $A \xcb B$ that $ \xCf A$ is
a proper subset of $B,$ likewise for $A \xcd B$ and $A \xcf B.$

Given some fixed set $U$ we work in, and $X \xcc U,$ then $ \xdC (X):=U-X$
.

If $ \xdy \xcc \xdp (X)$ for some
$X,$ we say that $ \xdy $ satisfies

$( \xcs )$ iff it is closed under finite intersections,

$( \xcS )$ iff it is closed under arbitrary intersections,

$( \xcv )$ iff it is closed under finite unions,

$( \xcV )$ iff it is closed under arbitrary unions,

$( \xdC )$ iff it is closed under complementation.

We will sometimes write $A=B \xFO C$ for: $A=B,$ or $A=C,$ or $A=B \xcv
C.$

We make ample and tacit use of the Axiom of Choice.
\index{Definition Rel-Base}

\ed

\bd

$\hspace{0.01em}$


\label{Definition Rel-Base}

$ \xeb^{*}$ will denote the transitive closure of the relation $ \xeb.$
If a relation $<,$
$ \xeb,$ or similar is given, $a \xcT b$ will express that a and $b$ are
$<-$ (or $ \xeb -)$
incomparable - context will tell. Given any relation $<,$ $ \xck $ will
stand for
$<$ or $=,$ conversely, given $ \xck,$ $<$ will stand for $ \xck,$ but
not $=,$ similarly
for $ \xeb $ etc.
\index{Definition Log-Base}

\ed

\bd

$\hspace{0.01em}$


\label{Definition Log-Base}

We work here in a classical propositional language $ \xdl,$ a theory $T$
will be an
arbitrary set of formulas. Formulas will often be named $ \xbf,$ $ \xbq
,$ etc., theories
$T,$ $S,$ etc.

$v( \xdl )$ will be the set of propositional variables of $ \xdl.$

$M_{ \xdl }$ will be the set of (classical) models of $ \xdl,$ $M(T)$ or
$M_{T}$
is the set of models of $T,$ likewise $M( \xbf )$ for a formula $ \xbf.$

$ \xdD_{ \xdl }:=\{M(T):$ $T$ a theory in $ \xdl \},$ the set of definable
model sets.

Note that, in classical propositional logic, $ \xCQ,M_{ \xdl } \xbe
\xdD_{ \xdl },$ $ \xdD_{ \xdl }$ contains
singletons, is closed under arbitrary intersections and finite unions.

An operation $f: \xdy \xcp \xdp (M_{ \xdl })$ for $ \xdy \xcc \xdp (M_{
\xdl })$ is called definability
preserving, $ \xCf (dp)$ or $( \xbm dp)$ in short, iff for all $X \xbe
\xdD_{ \xdl } \xcs \xdy $ $f(X) \xbe \xdD_{ \xdl }.$

We will also use $( \xbm dp)$ for binary functions $f: \xdy \xCK \xdy \xcp
\xdp (M_{ \xdl })$ - as needed
for theory revision - with the obvious meaning.

$ \xcl $ will be classical derivability, and

$ \ol{T}:=\{ \xbf:T \xcl \xbf \},$ the closure of $T$ under $ \xcl.$

$Con(.)$ will stand for classical consistency, so $Con( \xbf )$ will mean
that
$ \xbf $ is clasical consistent, likewise for $Con(T).$ $Con(T,T' )$ will
stand for
$Con(T \xcv T' ),$ etc.

Given a consequence relation $ \xcn,$ we define

$ \ol{ \ol{T} }:=\{ \xbf:T \xcn \xbf \}.$

(There is no fear of confusion with $ \ol{T},$ as it just is not useful to
close
twice under classical logic.)

$T \xco T':=\{ \xbf \xco \xbf ': \xbf \xbe T, \xbf ' \xbe T' \}.$

If $X \xcc M_{ \xdl },$ then $Th(X):=\{ \xbf:X \xcm \xbf \},$ likewise
for $Th(m),$ $m \xbe M_{ \xdl }.$
\index{Fact Log-Base}

\ed

We recollect and note:

\bfa

$\hspace{0.01em}$


\label{Fact Log-Base}

Let $ \xdl $ be a fixed propositional language, $ \xdD_{ \xdl } \xcc X,$ $
\xbm:X \xcp \xdp (M_{ \xdl }),$ for a $ \xdl -$theory $T$
$ \ol{ \ol{T} }:=Th( \xbm (M_{T})),$ let $T,$ $T' $ be arbitrary theories,
then:

(1) $ \xbm (M_{T}) \xcc M_{ \ol{ \ol{T} }}$,

(2) $M_{T} \xcv M_{T' }=M_{T \xco T' }$ and $M_{T \xcv T' }=M_{T} \xcs
M_{T' }$,

(3) $ \xbm (M_{T})= \xCQ $ $ \xcr $ $ \xcT \xbe \ol{ \ol{T} }$.

If $ \xbm $ is definability preserving or $ \xbm (M_{T})$ is finite, then
the following also hold:

(4) $ \xbm (M_{T})=M_{ \ol{ \ol{T} }}$,

(5) $T' \xcl \ol{ \ol{T} }$ $ \xcr $ $M_{T' } \xcc \xbm (M_{T}),$

(6) $ \xbm (M_{T})=M_{T' }$ $ \xcr $ $ \ol{T' }= \ol{ \ol{T} }.$
$ \xcz $
\\[3ex]
\index{Definition Log-Cond}

\efa

\bd

$\hspace{0.01em}$


\label{Definition Log-Cond}

We introduce here formally a list of properties of set functions on the
algebraic side, and their corresponding logical rules on the other side.

Recall that $ \ol{T}:=\{ \xbf:T \xcl \xbf \},$ $ \ol{ \ol{T} }:=\{ \xbf
:T \xcn \xbf \},$
where $ \xcl $ is classical consequence, and $ \xcn $ any other
consequence.

We show, wherever adequate, in parallel the formula version
in the left column, the theory version
in the middle column, and the semantical or algebraic
counterpart in the
right column. The algebraic counterpart gives conditions for a
function $f:\xdy\xcp\xdp (U)$, where $U$ is some set, and
$\xdy\xcc\xdp (U)$.

When the formula version is not commonly used, we omit it,
as we normally work only with the theory version.

Intuitively, $A$ and $B$ in the right hand side column stand for
$M(\xbf)$ for some formula $\xbf$, whereas $X$, $Y$ stand for
$M(T)$ for some theory $T$.

{\footnotesize

\begin{tabular}{|c|c|c|}

\hline

\multicolumn{3}{|c|}{Basics} \xEP

\hline

$(AND)$
\xEH
$(AND)$
\xEH
Closure under
\xEP

$ \xbf \xcn \xbq,  \xbf \xcn \xbq '   \xch $
\xEH
$ T \xcn \xbq, T \xcn \xbq '   \xch $
\xEH
finite
\xEP

$ \xbf \xcn \xbq \xcu \xbq ' $
\xEH
$ T \xcn \xbq \xcu \xbq ' $
\xEH
intersection
\xEP

\hline

$(OR)$ \xEH $(OR)$ \xEH $( \xbm OR)$ \xEP

$ \xbf \xcn \xbq,  \xbf ' \xcn \xbq   \xch $ \xEH
$ \ol{\ol{T}} \xcs \ol{\ol{T'}} \xcc \ol{\ol{T \xco T'}} $ \xEH
$f(X \xcv Y) \xcc f(X) \xcv f(Y)$
\xEP

$ \xbf \xco \xbf ' \xcn \xbq $ \xEH
\xEH
\xEP

\hline

$(wOR)$
\xEH
$(wOR)$
\xEH
$( \xbm wOR)$
\xEP

$ \xbf \xcn \xbq,$ $ \xbf ' \xcl \xbq $ $ \xch $
\xEH
$ \ol{ \ol{T} } \xcs \ol{T' }$ $ \xcc $ $ \ol{ \ol{T \xco T' } }$
\xEH
$f(X \xcv Y) \xcc f(X) \xcv Y$
\xEP

$ \xbf \xco \xbf ' \xcn \xbq $
\xEH
\xEH
\xEP

\hline

$(disjOR)$
\xEH
$(disjOR)$
\xEH
$( \xbm disjOR)$
\xEP

$ \xbf \xcl \xCN \xbf ',$ $ \xbf \xcn \xbq,$
\xEH
$\xCN Con(T \xcv T') \xch$
\xEH
$X \xcs Y= \xCQ $ $ \xch $
\xEP

$ \xbf ' \xcn \xbq $ $ \xch $ $ \xbf \xco \xbf ' \xcn \xbq $
\xEH
$ \ol{ \ol{T} } \xcs \ol{ \ol{T' } } \xcc \ol{ \ol{T \xco T' } }$
\xEH
$f(X \xcv Y) \xcc f(X) \xcv f(Y)$
\xEP

\hline

$(LLE)$
\xEH
$(LLE)$
\xEH
\xEP

Left Logical Equivalence
\xEH
\xEH
\xEP

$ \xcl \xbf \xcr \xbf ',  \xbf \xcn \xbq   \xch $
\xEH
$ \ol{T}= \ol{T' }  \xch   \ol{\ol{T}} = \ol{\ol{T'}}$
\xEH
trivially true
\xEP

$ \xbf ' \xcn \xbq $ \xEH \xEH \xEP

\hline

$(RW)$ Right Weakening
\xEH
$(RW)$
\xEH
upward closure
\xEP

$ \xbf \xcn \xbq,  \xcl \xbq \xcp \xbq '   \xch $
\xEH
$ T \xcn \xbq,  \xcl \xbq \xcp \xbq '   \xch $
\xEH
\xEP

$ \xbf \xcn \xbq ' $
\xEH
$T \xcn \xbq ' $
\xEH
\xEP

\hline

$(CCL)$ Classical Closure \xEH $(CCL)$ \xEH \xEP

\xEH
$ \ol{ \ol{T} }$ is classically
\xEH
trivially true
\xEP

\xEH closed \xEH \xEP

\hline

$(SC)$ Supraclassicality \xEH $(SC)$ \xEH $( \xbm \xcc )$ \xEP

$ \xbf \xcl \xbq $ $ \xch $ $ \xbf \xcn \xbq $ \xEH $ \ol{T} \xcc \ol{
\ol{T} }$ \xEH $f(X) \xcc X$ \xEP

\cline{1-1}

$(REF)$ Reflexivity \xEH \xEH \xEP
$ \xbD,\xba \xcn \xba $ \xEH \xEH \xEP

\hline

$(CP)$ \xEH $(CP)$ \xEH $( \xbm \xCQ )$ \xEP

Consistency Preservation \xEH \xEH \xEP

$ \xbf \xcn \xcT $ $ \xch $ $ \xbf \xcl \xcT $ \xEH $T \xcn \xcT $ $ \xch
$ $T \xcl \xcT $ \xEH $f(X)= \xCQ $ $ \xch $ $X= \xCQ $ \xEP

\hline

\xEH
\xEH $( \xbm \xCQ fin)$
\xEP

\xEH
\xEH $X \xEd \xCQ $ $ \xch $ $f(X) \xEd \xCQ $
\xEP

\xEH
\xEH for finite $X$
\xEP

\hline

\xEH $(PR)$ \xEH $( \xbm PR)$ \xEP

$ \ol{ \ol{ \xbf \xcu \xbf ' } }$ $ \xcc $ $ \ol{ \ol{ \ol{ \xbf } } \xcv
\{ \xbf ' \}}$ \xEH
$ \ol{ \ol{T \xcv T' } }$ $ \xcc $ $ \ol{ \ol{ \ol{T} } \xcv T' }$ \xEH
$X \xcc Y$ $ \xch $
\xEP

\xEH \xEH $f(Y) \xcs X \xcc f(X)$
\xEP

\cline{3-3}

\xEH
\xEH
$(\xbm PR ')$
\xEP

\xEH
\xEH
$f(X) \xcs Y \xcc f(X \xcs Y)$
\xEP

\hline

$(CUT)$ \xEH $(CUT)$ \xEH $ (\xbm CUT) $ \xEP
$ \xbD \xcn \xba; \xbD, \xba \xcn \xbb \xch $ \xEH
$T \xcc \ol{T' } \xcc \ol{ \ol{T} }  \xch $ \xEH
$f(X) \xcc Y \xcc X  \xch $ \xEP
$ \xbD \xcn \xbb $ \xEH
$ \ol{ \ol{T'} } \xcc \ol{ \ol{T} }$ \xEH
$f(X) \xcc f(Y)$
\xEP

\hline

\end{tabular}

}

{\footnotesize

\begin{tabular}{|c|c|c|}

\hline

\multicolumn{3}{|c|}{Cumulativity} \xEP

\hline

$(CM)$ Cautious Monotony \xEH $(CM)$ \xEH $ (\xbm CM) $ \xEP

$ \xbf \xcn \xbq,  \xbf \xcn \xbq '   \xch $ \xEH
$T \xcc \ol{T' } \xcc \ol{ \ol{T} }  \xch $ \xEH
$f(X) \xcc Y \xcc X  \xch $
\xEP

$ \xbf \xcu \xbq \xcn \xbq ' $ \xEH
$ \ol{ \ol{T} } \xcc \ol{ \ol{T' } }$ \xEH
$f(Y) \xcc f(X)$
\xEP

\cline{1-1}

\cline{3-3}

or $(ResM)$ Restricted Monotony \xEH \xEH $(\xbm ResM)$ \xEP
$ \xbD \xcn \xba, \xbb \xch \xbD,\xba \xcn \xbb $ \xEH \xEH
$ f(X) \xcc A \xcs B \xch f(X \xcs A) \xcc B $ \xEP

\hline

$(CUM)$ Cumulativity \xEH $(CUM)$ \xEH $( \xbm CUM)$ \xEP

$ \xbf \xcn \xbq   \xch $ \xEH
$T \xcc \ol{T' } \xcc \ol{ \ol{T} }  \xch $ \xEH
$f(X) \xcc Y \xcc X  \xch $
\xEP

$( \xbf \xcn \xbq '   \xcj   \xbf \xcu \xbq \xcn \xbq ' )$ \xEH
$ \ol{ \ol{T} }= \ol{ \ol{T' } }$ \xEH
$f(Y)=f(X)$ \xEP

\hline

\xEH
$ (\xcc \xcd) $
\xEH
$ (\xbm \xcc \xcd) $
\xEP
\xEH
$T \xcc \ol{\ol{T'}}, T' \xcc \ol{\ol{T}} \xch $
\xEH
$ f(X) \xcc Y, f(Y) \xcc X \xch $
\xEP
\xEH
$ \ol{\ol{T'}} = \ol{\ol{T}}$
\xEH
$ f(X)=f(Y) $
\xEP

\hline

\multicolumn{3}{|c|}{Rationality} \xEP

\hline

$(RatM)$ Rational Monotony \xEH $(RatM)$ \xEH $( \xbm RatM)$ \xEP

$ \xbf \xcn \xbq,  \xbf \xcN \xCN \xbq '   \xch $ \xEH
$Con(T \xcv \ol{\ol{T'}})$, $T \xcl T'$ $ \xch $ \xEH
$X \xcc Y, X \xcs f(Y) \xEd \xCQ   \xch $
\xEP

$ \xbf \xcu \xbq ' \xcn \xbq $ \xEH
$ \ol{\ol{T}} \xcd \ol{\ol{\ol{T'}} \xcv T} $ \xEH
$f(X) \xcc f(Y) \xcs X$ \xEP

\hline

\xEH $(RatM=)$ \xEH $( \xbm =)$ \xEP

\xEH
$Con(T \xcv \ol{\ol{T'}})$, $T \xcl T'$ $ \xch $ \xEH
$X \xcc Y, X \xcs f(Y) \xEd \xCQ   \xch $
\xEP

\xEH
$ \ol{\ol{T}} = \ol{\ol{\ol{T'}} \xcv T} $ \xEH
$f(X) = f(Y) \xcs X$ \xEP

\hline

\xEH
$(Log=' )$
\xEH $( \xbm =' )$
\xEP

\xEH
$Con( \ol{ \ol{T' } } \xcv T)$ $ \xch $
\xEH $f(Y) \xcs X \xEd \xCQ $ $ \xch $
\xEP

\xEH
$ \ol{ \ol{T \xcv T' } }= \ol{ \ol{ \ol{T' } } \xcv T}$
\xEH $f(Y \xcs X)=f(Y) \xcs X$
\xEP

\hline

\xEH
$(Log \xFO )$
\xEH $( \xbm \xFO )$
\xEP

\xEH
$ \ol{ \ol{T \xco T' } }$ is one of
\xEH $f(X \xcv Y)$ is one of
\xEP

\xEH
$\ol{\ol{T}},$ or $\ol{\ol{T'}},$ or $\ol{\ol{T}} \xcs \ol{\ol{T'}}$ (by (CCL))
\xEH $f(X),$ $f(Y)$ or $f(X) \xcv f(Y)$
\xEP

\hline

\xEH
$(Log \xcv )$
\xEH $( \xbm \xcv )$
\xEP

\xEH
$Con( \ol{ \ol{T' } } \xcv T),$ $ \xCN Con( \ol{ \ol{T' } }
\xcv \ol{ \ol{T} })$ $ \xch $
\xEH $f(Y) \xcs (X-f(X)) \xEd \xCQ $ $ \xch $
\xEP

\xEH
$ \xCN Con( \ol{ \ol{T \xco T' } } \xcv T' )$
\xEH $f(X \xcv Y) \xcs Y= \xCQ$
\xEP

\hline

\xEH
$(Log \xcv ' )$
\xEH $( \xbm \xcv ' )$
\xEP

\xEH
$Con( \ol{ \ol{T' } } \xcv T),$ $ \xCN Con( \ol{ \ol{T' }
} \xcv \ol{ \ol{T} })$ $ \xch $
\xEH $f(Y) \xcs (X-f(X)) \xEd \xCQ $ $ \xch $
\xEP

\xEH
$ \ol{ \ol{T \xco T' } }= \ol{ \ol{T} }$
\xEH $f(X \xcv Y)=f(X)$
\xEP

\hline

\xEH
\xEH $( \xbm \xbe )$
\xEP

\xEH
\xEH $a \xbe X-f(X)$ $ \xch $
\xEP

\xEH
\xEH $ \xcE b \xbe X.a \xce f(\{a,b\})$
\xEP

\hline

\end{tabular}

}

$(PR)$ is also called infinite conditionalization - we choose the name for
its central role for preferential structures $(PR)$ or $( \xbm PR).$

The system of rules $(AND)$ $(OR)$ $(LLE)$ $(RW)$ $(SC)$ $(CP)$ $(CM)$ $(CUM)$
is also called system $P$ (for preferential), adding $(RatM)$ gives the system
$R$ (for rationality or rankedness).

Roughly: Smooth preferential structures generate logics satisfying system
$P$, ranked structures logics satisfying system $R$.

A logic satisfying $(REF)$, $(ResM)$, and $(CUT)$ is called a consequence
relation.

$(LLE)$ and$(CCL)$ will hold automatically, whenever we work with model sets.

$(AND)$ is obviously closely related to filters, and corresponds to closure
under finite intersections. $(RW)$ corresponds to upward closure of filters.

More precisely, validity of both depend on the definition, and the
direction we consider.

Given $f$ and $(\xbm \xcc )$, $f(X)\xcc X$ generates a pricipal filter:
$\{X'\xcc X:f(X)\xcc X'\}$, with
the definition: If $X=M(T)$, then $T\xcn \xbf$  iff $f(X)\xcc M(\xbf )$.
Validity of $(AND)$ and
$(RW)$ are then trivial.

Conversely, we can define for $X=M(T)$

$\xdx:=\{X'\xcc X: \xcE \xbf (X'=X\xcs M(\xbf )$ and $T\xcn \xbf )\}$.

$(AND)$ then makes $\xdx$  closed under
finite intersections, $(RW)$ makes $\xdx$  upward
closed. This is in the infinite case usually not yet a filter, as not all
subsets of $X$ need to be definable this way.
In this case, we complete $\xdx$  by
adding all $X''$ such that there is $X'\xcc X''\xcc X$, $X'\xbe\xdx$.

Alternatively, we can define

$\xdx:=\{X'\xcc X: \xcS\{X \xcs M(\xbf ): T\xcn \xbf \} \xcc X' \}$.

$(SC)$ corresponds to the choice of a subset.

$(CP)$ is somewhat delicate, as it presupposes that the chosen model set is
non-empty. This might fail in the presence of ever better choices, without
ideal ones; the problem is addressed by the limit versions.

$(PR)$ is an infinitary version of one half of the deduction theorem: Let $T$
stand for $\xbf$, $T'$ for $\xbq$, and $\xbf \xcu \xbq \xcn \xbs$,
so $\xbf \xcn \xbq \xcp \xbs$, but $(\xbq \xcp \xbs )\xcu \xbq \xcl \xbs$.

$(CUM)$ (whose most interesting half in our context is $(CM)$) may best be seen
as
normal use of lemmas: We have worked hard and found some lemmas. Now
we can take a rest, and come back again with our new lemmas. Adding them to the
axioms will neither add new theorems, nor prevent old ones to hold.

\index{Fact Mu-Base}

\ed

\bfa

$\hspace{0.01em}$


\label{Fact Mu-Base}

This table is to be read as follows: If the left hand side holds for some
function $f: \xdy \xcp \xdp (U),$ and the auxiliary properties noted in
the middle also
hold for $f$ or $ \xdy,$ then the right hand side will hold, too - and
conversely.

{\small

\begin{tabular}{|c|c|c|c|}

\hline

\multicolumn{4}{|c|}{Basics} \xEP

\hline

(1.1)
\xEH
$(\xbm PR)$
\xEH
$\xch$ $(\xcs)+(\xbm \xcc)$
\xEH
$(\xbm PR')$
\xEP

\cline{1-1}

\cline{3-3}

(1.2)
\xEH
\xEH
$\xci$
\xEH
\xEP

\hline

(2.1)
\xEH
$(\xbm PR)$
\xEH
$\xch$ $(\xbm \xcc)$
\xEH
$(\xbm OR)$
\xEP

\cline{1-1}

\cline{3-3}

(2.2)
\xEH
\xEH
$\xci$ $(\xbm \xcc)$ + closure
\xEH
\xEP

\xEH
\xEH
under set difference
\xEH
\xEP

\hline

(3)
\xEH
$(\xbm PR)$
\xEH
$\xch$
\xEH
$( \xbm CUT)$
\xEP

\hline

(4)
\xEH
$(\xbm \xcc )+(\xbm \xcc \xcd )+(\xbm CUM)+$
\xEH
$\xcH$
\xEH
$( \xbm PR)$
\xEP

\xEH
$(\xbm RatM)+(\xcs )$
\xEH
\xEH
\xEP

\hline

\multicolumn{4}{|c|}{Cumulativity} \xEP

\hline

(5.1)
\xEH
$(\xbm CM)$
\xEH
$\xch$ $(\xcs)+(\xbm \xcc)$
\xEH
$(\xbm ResM)$
\xEP

\cline{1-1}

\cline{3-3}

(5.2)
\xEH
\xEH
$\xci$ (infin.)
\xEH
\xEP

\hline

(6)
\xEH
$(\xbm CM)+(\xbm CUT)$
\xEH
$\xcj$
\xEH
$(\xbm CUM)$
\xEP

\hline

(7)
\xEH
$( \xbm \xcc )+( \xbm \xcc \xcd )$
\xEH
$\xch$
\xEH
$( \xbm CUM)$
\xEP

\hline

(8)
\xEH
$( \xbm \xcc )+( \xbm CUM)+( \xcs )$
\xEH
$\xch$
\xEH
$( \xbm \xcc \xcd )$
\xEP

\hline

(9)
\xEH
$( \xbm \xcc )+( \xbm CUM)$
\xEH
$\xcH$
\xEH
$( \xbm \xcc \xcd )$
\xEP

\hline

\multicolumn{4}{|c|}{Rationality} \xEP

\hline

(10)
\xEH
$( \xbm RatM )+( \xbm PR )$
\xEH
$\xch$
\xEH
$( \xbm =)$
\xEP

\hline

(11)
\xEH
$( \xbm =)$
\xEH
$ \xch $
\xEH
$( \xbm PR),$
\xEP

\hline

(12.1)
\xEH
$( \xbm =)$
\xEH
$ \xch $ $(\xcs)+( \xbm \xcc )$
\xEH
$( \xbm =' ),$
\xEP
\cline{1-1}
\cline{3-3}
(12.2)
\xEH
\xEH
$ \xci $
\xEH
\xEP

\hline

(13)
\xEH
$( \xbm \xcc ),$ $( \xbm =)$
\xEH
$ \xch $ $(\xcv)$
\xEH
$( \xbm \xcv ),$
\xEP

\hline

(14)
\xEH
$( \xbm \xcc ),$ $( \xbm \xCQ ),$ $( \xbm =)$
\xEH
$ \xch $ $(\xcv)$
\xEH
$( \xbm \xFO ),$ $( \xbm \xcv ' ),$ $( \xbm CUM),$
\xEP

\hline

(15)
\xEH
$( \xbm \xcc )+( \xbm \xFO )$
\xEH
$ \xch $ $\xdy$ closed under set difference
\xEH
$( \xbm =),$
\xEP

\hline

(16)
\xEH
$( \xbm \xFO )+( \xbm \xbe )+( \xbm PR)+$
\xEH
$ \xch $ $(\xcv)$ + $\xdy$ contains singletons
\xEH
$( \xbm =),$
\xEP
\xEH
$( \xbm \xcc )$
\xEH
\xEH
\xEP

\hline

(17)
\xEH
$( \xbm CUM)+( \xbm =)$
\xEH
$ \xch $ $(\xcv)$ + $\xdy$ contains singletons
\xEH
$( \xbm \xbe ),$
\xEP

\hline

(18)
\xEH
$( \xbm CUM)+( \xbm =)+( \xbm \xcc )$
\xEH
$ \xch $ $(\xcv)$
\xEH
$( \xbm \xFO ),$
\xEP

\hline

(19)
\xEH
$( \xbm PR)+( \xbm CUM)+( \xbm \xFO )$
\xEH
$ \xch $ sufficient, e.g. true in $\xdD_{\xdl}$
\xEH
$( \xbm =)$.
\xEP

\hline

(20)
\xEH
$( \xbm \xcc )+( \xbm PR)+( \xbm =)$
\xEH
$ \xcH $
\xEH
$( \xbm \xFO ),$
\xEP

\hline

(21)
\xEH
$( \xbm \xcc )+( \xbm PR)+( \xbm \xFO )$
\xEH
$ \xcH $ (without closure
\xEH
$( \xbm =)$
\xEP
\xEH
\xEH
under set difference),
\xEH
\xEP

\hline

(22)
\xEH
$( \xbm \xcc )+( \xbm PR)+( \xbm \xFO )+$
\xEH
$ \xcH $
\xEH
$( \xbm \xbe )$
\xEP
\xEH
$( \xbm =)+( \xbm \xcv )$
\xEH
\xEH
(thus not representability
\xEP
\xEH
\xEH
\xEH
by ranked structures)
\xEP

\hline

\end{tabular}

}

\index{Proposition Alg-Log}

\efa

\bp

$\hspace{0.01em}$


\label{Proposition Alg-Log}

The following table is to be read as follows:

Let a logic $ \xcn $ satisfies $ \xCf (LLE)$ and $ \xCf (CCL),$ and define
a function $f: \xdD_{ \xdl } \xcp \xdD_{ \xdl }$
by $f(M(T)):=M( \ol{ \ol{T} }).$ Then $f$ is well defined, satisfies $(
\xbm dp),$ and $ \ol{ \ol{T} }=Th(f(M(T))).$

If $ \xcn $ satisfies a rule in the left hand side,
then - provided the additional properties noted in the middle for $ \xch $
hold, too -
$f$ will satisfy the property in the right hand side.

Conversely, if $f: \xdy \xcp \xdp (M_{ \xdl })$ is a function, with $
\xdD_{ \xdl } \xcc \xdy,$ and we define a logic
$ \xcn $ by $ \ol{ \ol{T} }:=Th(f(M(T))),$ then $ \xcn $ satisfies $ \xCf
(LLE)$ and $ \xCf (CCL).$
If $f$ satisfies $( \xbm dp),$ then $f(M(T))=M( \ol{ \ol{T} }).$

If $f$ satisfies a property in the right hand side,
then - provided the additional properties noted in the middle for $ \xci $
hold, too -
$ \xcn $ will satisfy the property in the left hand side.

If ``formula'' is noted in the table, this means that, if one of the
theories
(the one named the same way in Definition \ref{Definition Log-Cond})
is equivalent to a formula, we can renounce on $( \xbm dp).$

{\small

\begin{tabular}{|c|c|c|c|}

\hline

\multicolumn{4}{|c|}{Basics} \xEP

\hline

(1.1) \xEH $(OR)$ \xEH $\xch$ \xEH $(\xbm OR)$ \xEP

\cline{1-1}

\cline{3-3}

(1.2) \xEH \xEH $\xci$ \xEH \xEP

\hline

(2.1) \xEH $(disjOR)$ \xEH $\xch$ \xEH $(\xbm disjOR)$ \xEP

\cline{1-1}

\cline{3-3}

(2.2) \xEH \xEH $\xci$ \xEH \xEP

\hline

(3.1) \xEH $(wOR)$ \xEH $\xch$ \xEH $(\xbm wOR)$ \xEP

\cline{1-1}

\cline{3-3}

(3.2) \xEH \xEH $\xci$ \xEH \xEP

\hline

(4.1) \xEH $(SC)$ \xEH $\xch$ \xEH $(\xbm \xcc)$ \xEP

\cline{1-1}

\cline{3-3}

(4.2) \xEH \xEH $\xci$ \xEH \xEP

\hline

(5.1) \xEH $(CP)$ \xEH $\xch$ \xEH $(\xbm \xCQ)$ \xEP

\cline{1-1}

\cline{3-3}

(5.2) \xEH \xEH $\xci$ \xEH \xEP

\hline

(6.1) \xEH $(PR)$ \xEH $\xch$ \xEH $(\xbm PR)$ \xEP

\cline{1-1}

\cline{3-3}

(6.2) \xEH \xEH $\xci$ $(\xbm dp)+(\xbm \xcc)$ \xEH \xEP

\cline{1-1}

\cline{3-3}

(6.3) \xEH \xEH $\xcI$ without $(\xbm dp)$ \xEH \xEP

\cline{1-1}

\cline{3-3}

(6.4) \xEH \xEH $\xci$ $(\xbm \xcc)$ \xEH \xEP

\xEH \xEH $T'$ a formula \xEH \xEP

\hline

(6.5) \xEH $(PR)$ \xEH $\xci$ \xEH $(\xbm PR')$ \xEP

\xEH \xEH $T'$ a formula \xEH \xEP

\hline

(7.1) \xEH $(CUT)$ \xEH $\xch$ \xEH $(\xbm CUT)$ \xEP

\cline{1-1}

\cline{3-3}

(7.2) \xEH \xEH $\xci$ \xEH \xEP

\hline

\multicolumn{4}{|c|}{Cumulativity} \xEP

\hline

(8.1) \xEH $(CM)$ \xEH $\xch$ \xEH $(\xbm CM)$ \xEP

\cline{1-1}

\cline{3-3}

(8.2) \xEH \xEH $\xci$ \xEH \xEP

\hline

(9.1) \xEH $(ResM)$ \xEH $\xch$ \xEH $(\xbm ResM)$ \xEP

\cline{1-1}

\cline{3-3}

(9.2) \xEH \xEH $\xci$ \xEH \xEP

\hline

(10.1) \xEH $(\xcc \xcd)$ \xEH $\xch$ \xEH $(\xbm \xcc \xcd)$ \xEP

\cline{1-1}

\cline{3-3}

(10.2) \xEH \xEH $\xci$ \xEH \xEP

\hline

(11.1) \xEH $(CUM)$ \xEH $\xch$ \xEH $(\xbm CUM)$ \xEP

\cline{1-1}

\cline{3-3}

(11.2) \xEH \xEH $\xci$ \xEH \xEP

\hline

\multicolumn{4}{|c|}{Rationality} \xEP

\hline

(12.1) \xEH $(RatM)$ \xEH $\xch$ \xEH $(\xbm RatM)$ \xEP

\cline{1-1}

\cline{3-3}

(12.2) \xEH \xEH $\xci$ $(\xbm dp)$ \xEH \xEP

\cline{1-1}

\cline{3-3}

(12.3) \xEH \xEH $\xcI$ without $(\xbm dp)$ \xEH \xEP

\cline{1-1}

\cline{3-3}

(12.4) \xEH \xEH $\xci$ \xEH \xEP

\xEH \xEH $T$ a formula \xEH \xEP

\hline

(13.1) \xEH $(RatM=)$ \xEH $\xch$ \xEH $(\xbm =)$ \xEP

\cline{1-1}

\cline{3-3}

(13.2) \xEH \xEH $\xci$ $(\xbm dp)$ \xEH \xEP

\cline{1-1}

\cline{3-3}

(13.3) \xEH \xEH $\xcI$ without $(\xbm dp)$ \xEH \xEP

\cline{1-1}

\cline{3-3}

(13.4) \xEH \xEH $\xci$ \xEH \xEP

\xEH \xEH $T$ a formula \xEH \xEP

\hline

(14.1) \xEH $(Log = ')$ \xEH $\xch$ \xEH $(\xbm = ')$ \xEP

\cline{1-1}
\cline{3-3}

(14.2) \xEH \xEH $\xci$ $(\xbm dp)$ \xEH \xEP

\cline{1-1}
\cline{3-3}

(14.3) \xEH \xEH $\xcI$ without $(\xbm dp)$ \xEH \xEP

\cline{1-1}
\cline{3-3}

(14.4) \xEH \xEH $\xci$ $T$ a formula \xEH \xEP

\hline

(15.1) \xEH $(Log \xFO )$ \xEH $\xch$ \xEH $(\xbm \xFO )$ \xEP

\cline{1-1}
\cline{3-3}

(15.2) \xEH \xEH $\xci$ \xEH \xEP

\hline

(16.1)
\xEH
$(Log \xcv )$
\xEH
$\xch$ $(\xbm \xcc)+(\xbm =)$
\xEH
$(\xbm \xcv )$
\xEP

\cline{1-1}
\cline{3-3}

(16.2) \xEH \xEH $\xci$ $(\xbm dp)$ \xEH \xEP

\cline{1-1}
\cline{3-3}

(16.3) \xEH \xEH $\xcI$ without $(\xbm dp)$ \xEH \xEP

\hline

(17.1)
\xEH
$(Log \xcv ')$
\xEH
$\xch$ $(\xbm \xcc)+(\xbm =)$
\xEH
$(\xbm \xcv ')$
\xEP

\cline{1-1}
\cline{3-3}

(17.2) \xEH \xEH $\xci$ $(\xbm dp)$ \xEH \xEP

\cline{1-1}
\cline{3-3}

(17.3) \xEH \xEH $\xcI$ without $(\xbm dp)$ \xEH \xEP

\hline

\end{tabular}

}

\index{Definition Pref-Str}

\ep

\bd

$\hspace{0.01em}$


\label{Definition Pref-Str}

Fix $U \xEd \xCQ,$ and consider arbitrary $X.$
Note that this $X$ has not necessarily anything to do with $U,$ or $ \xdu
$ below.
Thus, the functions $ \xbm_{ \xdm }$ below are in principle functions from
$V$ to $V$ - where $V$
is the set theoretical universe we work in.

(A) Preferential models or structures.

(1) The version without copies:

A pair $ \xdm:=<U, \xeb >$ with $U$ an arbitrary set, and $ \xeb $ an
arbitrary binary relation
is called a preferential model or structure.

(2) The version with copies:

A pair $ \xdm:=< \xdu, \xeb >$ with $ \xdu $ an arbitrary set of pairs,
and $ \xeb $ an arbitrary binary
relation is called a preferential model or structure.

If $<x,i> \xbe \xdu,$ then $x$ is intended to be an element of $U,$ and
$i$ the index of the
copy.

We sometimes also need copies of the relation $ \xeb,$ we will then
replace $ \xeb $
by one or several arrows $ \xba $ attacking non-minimal elements, e.g. $x
\xeb y$ will
be written $ \xba:x \xcp y,$ $<x,i> \xeb <y,i>$ will be written $ \xba
:<x,i> \xcp <y,i>,$ and
finally we might have $< \xba,k>:x \xcp y$ and $< \xba,k>:<x,i> \xcp
<y,i>,$ etc.

(B) Minimal elements, the functions $ \xbm_{ \xdm }$

(1) The version without copies:

Let $ \xdm:=<U, \xeb >,$ and define

$ \xbm_{ \xdm }(X)$ $:=$ $\{x \xbe X:$ $x \xbe U$ $ \xcu $ $ \xCN \xcE x'
\xbe X \xcs U.x' \xeb x\}.$

$ \xbm_{ \xdm }(X)$ is called the set of minimal elements of $X$ (in $
\xdm ).$

(2) The version with copies:

Let $ \xdm:=< \xdu, \xeb >$ be as above. Define

$ \xbm_{ \xdm }(X)$ $:=$ $\{x \xbe X:$ $ \xcE <x,i> \xbe \xdu. \xCN \xcE
<x',i' > \xbe \xdu (x' \xbe X$ $ \xcu $ $<x',i' >' \xeb <x,i>)\}.$

Again, by abuse of language, we say that $ \xbm_{ \xdm }(X)$ is the set of
minimal elements
of $X$ in the structure. If the context is clear, we will also write just
$ \xbm.$

We sometimes say that $<x,i>$ ``kills'' or ``minimizes'' $<y,j>$ if
$<x,i> \xeb <y,j>.$ By abuse of language we also say a set $X$ kills or
minimizes a set
$Y$ if for all $<y,j> \xbe \xdu,$ $y \xbe Y$ there is $<x,i> \xbe \xdu,$
$x \xbe X$ s.t. $<x,i> \xeb <y,j>.$

$ \xdm $ is also called injective or 1-copy, iff there is always at most
one copy
$<x,i>$ for each $x.$ Note that the existence of copies corresponds to a
non-injective labelling function - as is often used in nonclassical
logic, e.g. modal logic.

We say that $ \xdm $ is transitive, irreflexive, etc., iff $ \xeb $ is.

Note that $ \xbm (X)$ might well be empty, even if $X$ is not.
\index{Definition Pref-Log}

\ed

\bd

$\hspace{0.01em}$


\label{Definition Pref-Log}

We define the consequence relation of a preferential structure for a
given propositional language $ \xdl.$

(A)

(1) If $m$ is a classical model of a language $ \xdl,$ we say by abuse
of language

$<m,i> \xcm \xbf $ iff $m \xcm \xbf,$

and if $X$ is a set of such pairs, that

$X \xcm \xbf $ iff for all $<m,i> \xbe X$ $m \xcm \xbf.$

(2) If $ \xdm $ is a preferential structure, and $X$ is a set of $ \xdl
-$models for a
classical propositional language $ \xdl,$ or a set of pairs $<m,i>,$
where the $m$ are
such models, we call $ \xdm $ a classical preferential structure or model.

(B)

Validity in a preferential structure, or the semantical consequence
relation
defined by such a structure:

Let $ \xdm $ be as above.

We define:

$T \xcm_{ \xdm } \xbf $ iff $ \xbm_{ \xdm }(M(T)) \xcm \xbf,$ i.e. $
\xbm_{ \xdm }(M(T)) \xcc M( \xbf ).$

$ \xdm $ will be called definability preserving iff for all $X \xbe \xdD_{
\xdl }$ $ \xbm_{ \xdm }(X) \xbe \xdD_{ \xdl }.$

As $ \xbm_{ \xdm }$ is defined on $ \xdD_{ \xdl },$ but need by no means
always result in some new
definable set, this is (and reveals itself as a quite strong) additional
property.
\index{Definition Smooth}

\ed

\bd

$\hspace{0.01em}$


\label{Definition Smooth}

Let $ \xdy \xcc \xdp (U).$ (In applications to logic, $ \xdy $ will be $
\xdD_{ \xdl }.)$

A preferential structure $ \xdm $ is called $ \xdy -$smooth iff in every
$X \xbe \xdy $ every element
$x \xbe X$ is either minimal in $X$ or above an element, which is minimal
in $X.$ More
precisely:

(1) The version without copies:

If $x \xbe X \xbe \xdy,$ then either $x \xbe \xbm (X)$ or there is $x'
\xbe \xbm (X).x' \xeb x.$

(2) The version with copies:

If $x \xbe X \xbe \xdy,$ and $<x,i> \xbe \xdu,$ then either there is no
$<x',i' > \xbe \xdu,$ $x' \xbe X,$
$<x',i' > \xeb <x,i>$ or there is $<x',i' > \xbe \xdu,$ $<x',i' > \xeb
<x,i>,$ $x' \xbe X,$ s.t. there is
no $<x'',i'' > \xbe \xdu,$ $x'' \xbe X,$ with $<x'',i'' > \xeb <x',i'
>.$

When considering the models of a language $ \xdl,$ $ \xdm $ will be
called smooth iff
it is $ \xdD_{ \xdl }-$smooth; $ \xdD_{ \xdl }$ is the default.

Obviously, the richer the set $ \xdy $ is, the stronger the condition $
\xdy -$smoothness
will be.
\index{Table Pref-Representation-Without-Ref}

\ed

The following table summarizes representation by not necessarily ranked
preferential structures. The implications on the right are shown in
Proposition \ref{Proposition Alg-Log} (going via the $ \xbm -$functions),
those on the left
are shown in the respective representation theorems.
\label{Table Pref-Representation-Without-Ref}

{\scriptsize

\begin{tabular}{|c|c|c|c|c|}

\hline

$\xbm-$ function
\xEH
\xEH
Pref.Structure
\xEH
\xEH
Logic
\xEP

\hline

$(\xbm \xcc)+(\xbm PR)$
\xEH
$\xci$
\xEH
general
\xEH
$\xch$ $(\xbm dp)$
\xEH
$(LLE)+(RW)+(SC)+(PR)$
\xEP

\cline{2-2}
\cline{4-4}

\xEH
$\xch$
\xEH
\xEH
$\xci$
\xEH
\xEP

\cline{2-2}
\cline{4-4}

\xEH
\xEH
\xEH
$\xcH$ without $(\xbm dp)$
\xEH
\xEP

\hline

$(\xbm \xcc)+(\xbm PR)$
\xEH
$\xci$
\xEH
transitive
\xEH
$\xch$ $(\xbm dp)$
\xEH
$(LLE)+(RW)+(SC)+(PR)$
\xEP

\cline{2-2}
\cline{4-4}

\xEH
$\xch$
\xEH
\xEH
$\xci$
\xEH
\xEP

\cline{2-2}
\cline{4-4}

\xEH
\xEH
\xEH
$\xcH$ without $(\xbm dp)$
\xEH
\xEP

\hline

$(\xbm \xcc)+(\xbm PR)+(\xbm CUM)$
\xEH
$\xci$
\xEH
smooth
\xEH
$\xch$ $(\xbm dp)$
\xEH
$(LLE)+(RW)+(SC)+(PR)+$
\xEP

\xEH
\xEH
\xEH
\xEH
$(CUM)$
\xEP

\cline{2-2}
\cline{4-4}

\xEH
$\xch$ $(\xcv)$
\xEH
\xEH
$\xci$ $(\xcv)$
\xEH
\xEP

\cline{2-2}
\cline{4-4}

\xEH
\xEH
\xEH
$\xcH$ without $(\xbm dp)$
\xEH
\xEP

\hline

$(\xbm \xcc)+(\xbm PR)+(\xbm CUM)$
\xEH
$\xci$
\xEH
smooth+transitive
\xEH
$\xch$ $(\xbm dp)$
\xEH
$(LLE)+(RW)+(SC)+(PR)+$
\xEP

\xEH
\xEH
\xEH
\xEH
$(CUM)$
\xEP

\cline{2-2}
\cline{4-4}

\xEH
$\xch$ $(\xcv)$
\xEH
\xEH
$\xci$ $(\xcv)$
\xEH
\xEP

\cline{2-2}
\cline{4-4}

\xEH
\xEH
\xEH
$\xcH$ without $(\xbm dp)$
\xEH
\xEP

\hline

\end{tabular}

}

\subsection{
Generalized preferential structures
}


\label{Section Reac-GenPref}
\index{Section Reac-GenPref}

\subsubsection{
Introduction
}


\label{Section Reac-GenPref-Intro}
\index{Section Reac-GenPref-Intro}
\index{Comment Gen-Pref}

\bcom

$\hspace{0.01em}$


\label{Comment Gen-Pref}

A counterargument to $ \xba $ is NOT an argument for $ \xCN \xba $ (this
is asking for too
much), but just showing one case where $ \xCN \xba $ holds. In
preferential structures,
an argument for $ \xba $ is a set of level 1 arrows, eliminating $ \xCN
\xba -$models. A
counterargument is one level 2 arrow, attacking one such level 1 arrow.

Of course, when we have copies, we may need many successful attacks, on
all
copies, to achieve the goal. As we may have copies of level 1 arrows, we
may need many level 2 arrows to destroy them all.
\index{Definition Generalized preferential structure}

\ecom

\bd

$\hspace{0.01em}$


\label{Definition Generalized preferential structure}

An IBR is called a generalized preferential structure iff the origins of
all
arrows are points. We will usually write $x,y$ etc. for points, $ \xba,$
$ \xbb $ etc. for
arrows.
\index{Definition Level-n-Arrow}

\ed

\bd

$\hspace{0.01em}$


\label{Definition Level-n-Arrow}

Consider a generalized preferential structure $ \xdx.$

(1) Level $n$ arrow:

Definition by upward induction.

If $ \xba:x \xcp y,$ $x,y$ are points, then $ \xba $ is a level 1 arrow.

If $ \xba:x \xcp \xbb,$ $x$ is a point, $ \xbb $ a level $n$ arrow, then
$ \xba $ is a level $n+1$ arrow.
$(o( \xba )$ is the origin, $d( \xba )$ is the destination of $ \xba.)$

$ \xbl ( \xba )$ will denote the level of $ \xba.$

(2) Level $n$ structure:

$ \xdx $ is a level $n$ structure iff all arrows in $ \xdx $ are at most
level $n$ arrows.

We consider here only structures of some arbitrary but finite level $n.$

(3) We define for an arrow $ \xba $ by induction $O( \xba )$ and $D( \xba
).$

If $ \xbl ( \xba )=1,$ then $O( \xba ):=\{o( \xba )\},$ $D( \xba ):=\{d(
\xba )\}.$

If $ \xba:x \xcp \xbb,$ then $D( \xba ):=D( \xbb ),$ and $O( \xba
):=\{x\} \xcv O( \xbb ).$

Thus, for example, if $ \xba:x \xcp y,$ $ \xbb:z \xcp \xba,$ then $O(
\xbb ):=\{x,z\},$ $D( \xbb )=\{y\}.$
\index{Example Inf-Level}

\ed

We will not consider here diagrams with arbitrarily high levels. One
reason is
that diagrams like the following will have an unclear meaning:

\be

$\hspace{0.01em}$


\label{Example Inf-Level}

$< \xba,1>:x \xcp y,$

$< \xba,n+1>:x \xcp < \xba,n>$ $(n \xbe \xbo ).$

Is $y \xbe \xbm (X)?$
\index{Definition Valid-Arrow}

\ee

\bd

$\hspace{0.01em}$


\label{Definition Valid-Arrow}

Let $ \xdx $ be a generalized preferential structure of (finite) level
$n.$

We define (by downward induction):

(1) Valid $X-to-Y$ arrow:

Let $X,Y \xcc \xdP ( \xdx ).$

$ \xba \xbe \xdA ( \xdx )$ is a valid $X-to-Y$ arrow iff

(1.1) $O( \xba ) \xcc X,$ $D( \xba ) \xcc Y,$

(1.2) $ \xcA \xbb:x' \xcp \xba.(x' \xbe X$ $ \xch $ $ \xcE \xbg:x''
\xcp \xbb.( \xbg $ is a valid $X-to-Y$ arrow)).

We will also say that $ \xba $ is a valid arrow in $X,$ or just valid in
$X,$ iff $ \xba $ is a
valid $X-to-X$ arrow.

(2) Valid $X \xch Y$ arrow:

Let $X \xcc Y \xcc \xdP ( \xdx ).$

$ \xba \xbe \xdA ( \xdx )$ is a valid $X \xch Y$ arrow iff

(2.1) $o( \xba ) \xbe X,$ $O( \xba ) \xcc Y,$ $D( \xba ) \xcc Y,$

(2.2) $ \xcA \xbb:x' \xcp \xba.(x' \xbe Y$ $ \xch $ $ \xcE \xbg:x''
\xcp \xbb.( \xbg $ is a valid $X \xch Y$ arrow)).

(Note that in particular $o( \xbg ) \xbe X,$ and that $o( \xbb )$ need not
be in $X,$ but
can be in the bigger $Y.)$
\index{Fact Higher-Validity}

\ed

\bfa

$\hspace{0.01em}$


\label{Fact Higher-Validity}

(1) If $ \xba $ is a valid $X \xch Y$ arrow, then $ \xba $ is a valid
$Y-to-Y$ arrow.

(2) If $X \xcc X' \xcc Y' \xcc Y \xcc \xdP ( \xdx )$ and $ \xba \xbe \xdA
( \xdx )$ is a valid $X \xch Y$ arrow, and
$O( \xba ) \xcc Y',$ $D( \xba ) \xcc Y',$ then $ \xba $ is a valid $X'
\xch Y' $ arrow.
\index{Fact Higher-Validity Proof}

\efa

\paragraph{
Proof Fact Higher-Validity
}

$\hspace{0.01em}$


\label{Section Proof Fact Higher-Validity}

Let $ \xba $ be a valid $X \xch Y$ arrow. We show (1) and (2) together by
downward induction
(both are trivial).

By prerequisite $o( \xba ) \xbe X \xcc X',$ $O( \xba ) \xcc Y' \xcc Y,$
$D( \xba ) \xcc Y' \xcc Y.$

Case 1: $ \xbl ( \xba )=n.$ So $ \xba $ is a valid $X' \xch Y' $ arrow,
and a valid $Y-to-Y$ arrow.

Case 2: $ \xbl ( \xba )=n-1.$ So there is no $ \xbb:x' \xcp \xba,$ $y
\xbe Y,$ so $ \xba $ is a valid
$Y-to-Y$ arrow. By $Y' \xcc Y$ $ \xba $ is a valid $X' \xch Y' $ arrow.

Case 3: Let the result be shown down to $m,$ $n>m>1,$ let $ \xbl ( \xba
)=m-1.$
So $ \xcA \xbb:x' \xcp \xba (x' \xbe Y$ $ \xch $ $ \xcE \xbg:x'' \xcp
\xbb (x'' \xbe X$ and $ \xbg $ is a valid $X \xch Y$ arrow)).
By induction hypothesis $ \xbg $ is a valid $Y-to-Y$ arrow, and a valid
$X' \xch Y' $ arrow. So $ \xba $ is a valid $Y-to-Y$ arrow, and by $Y'
\xcc Y,$ $ \xba $ is a valid
$X' \xch Y' $ arrow.

$ \xcz $
\\[3ex]
\index{Definition Higher-Mu}

\bd

$\hspace{0.01em}$


\label{Definition Higher-Mu}

Let $ \xdx $ be a generalized preferential structure of level $n,$ $X \xcc
\xdP ( \xdx ).$

$ \xbm (X):=\{x \xbe X:$ $ \xcE <x,i>. \xCN \xcE $ valid $X-to-X$ arrow $
\xba:x' \xcp <x,i>\}.$
\index{Comment Smooth-Gen}

\ed

\bcom

$\hspace{0.01em}$


\label{Comment Smooth-Gen}

The purpose of smoothness is to guarantee cumulativity. Smoothness
achieves Cumulativity by mirroring all information present in $X$ also in
$ \xbm (X).$
Closer inspection shows that smoothness does more than necessary. This is
visible when there are copies (or, equivalently, non-injective labelling
functions). Suppose we have two copies of $x \xbe X,$ $<x,i>$ and $<x,i'
>,$ and there
is $y \xbe X,$ $ \xba:<y,j> \xcp <x,i>,$ but there is no $ \xba ':<y'
,j' > \xcp <x,i' >,$ $y' \xbe X.$ Then
$ \xba:<y,j> \xcp <x,i>$
is irrelevant, as $x \xbe \xbm (X)$ anyhow. So mirroring $ \xba:<y,j>
\xcp <x,i>$ in $ \xbm (X)$ is not
necessary, i.e. it is not necessary to have some $ \xba ':<y',j' > \xcp
<x,i>,$ $y' \xbe \xbm (X).$

On the other hand, Example \ref{Example Need-Smooth} shows that,
if we want smooth structures to correspond to the property $( \xbm CUM),$
we
need at least some valid arrows from $ \xbm (X)$ also for higher level
arrows.
This ``some'' is made precise (essentially) in Definition \ref{Definition
X-Sub-X'}.

From a more philosophical point of view,
when we see the (inverted) arrows of preferential structures as attacks on
non-minimal elements, then we should see smooth structures as always
having
attacks also from valid (minimal) elements. So, in general structures,
also
attacks from non-valid elements are valid, in smooth structures we always
also have attacks from valid elements.

The analogon to usual smooth structures, on level 2, is then that any
successfully attacked level 1 arrow is also attacked from a minimal point.
\index{Definition X-Sub-X'}

\ecom

\bd

$\hspace{0.01em}$


\label{Definition X-Sub-X'}

Let $ \xdx $ be a generalized preferential structure.

$X \xes X' $ iff

(1) $X \xcc X' \xcc \xdP ( \xdx ),$

(2) $ \xcA x \xbe X' -X$ $ \xcA <x,i>$ $ \xcE \xba:x' \xcp <x,i>( \xba $
is a valid $X \xch X' $ arrow),

(3) $ \xcA x \xbe X$ $ \xcE <x,i>$

$ \xDC $ $( \xcA \xba:x' \xcp <x,i>(x' \xbe X' $ $ \xch $ $ \xcE \xbb
:x'' \xcp \xba.( \xbb $ is a valid $X \xch X' $ arrow))).

Note that (3) is not simply the negation of (2):

Consider a level 1 structure. Thus all level 1 arrows are valid, but the
source
of the arrows must not be neglected.

(2) reads now: $ \xcA x \xbe X' -X$ $ \xcA <x,i>$ $ \xcE \xba:x' \xcp
<x,i>.x' \xbe X$

(3) reads: $ \xcA x \xbe X$ $ \xcE <x,i>$ $ \xCN \xcE \xba:x' \xcp
<x,i>.x' \xbe X' $

This is intended: intuitively, $X= \xbm (X' ),$ and minimal elements must
not be
attacked at all, but non-minimals must be attacked from $X$ - which is a
modified
version of smoothness.
\index{Remark X-Sub-X'}

\ed

\br

$\hspace{0.01em}$


\label{Remark X-Sub-X'}

We note the special case of Definition \ref{Definition X-Sub-X'} for level
3 structures,
as it will be used later. We also write it immediately for the intended
case
$ \xbm (X) \xes X,$ and explicitly with copies.

$x \xbe \xbm (X)$ iff

(1) $ \xcE <x,i> \xcA < \xba,k>:<y,j> \xcp <x,i>$

$ \xDC (y \xbe X$ $ \xcp $ $ \xcE < \xbb ',l' >:<z',m' > \xcp < \xba
,k>.$

$ \xDC \xDC (z' \xbe \xbm (X)$ $ \xcu $ $ \xCN \xcE < \xbg ',n' >:<u',p'
> \xcp < \xbb ',l' >.u' \xbe X))$

See Diagram \ref{Diagram Essential-Smooth-3-1-2}.

$x \xbe X- \xbm (X)$ iff

(2) $ \xcA <x,i> \xcE < \xba ',k' >:<y',j' > \xcp <x,i>$

$ \xDC (y' \xbe \xbm (X)$ $ \xcu $

$ \xDC \xDC (a)$ $ \xCN \xcE < \xbb ',l' >:<z',m' > \xcp < \xba ',k'
>.z' \xbe X$

$ \xDC \xDC or$

$ \xDC \xDC (b)$ $ \xcA < \xbb ',l' >:<z',m' > \xcp < \xba ',k' >$

$ \xDC \xDC \xDC (z' \xbe X$ $ \xcp $ $ \xcE < \xbg ',n' >:<u',p' > \xcp
< \xbb ',l' >.u' \xbe \xbm (X))$ )

See Diagram \ref{Diagram Essential-Smooth-3-2}.

\vspace{15mm}

\begin{diagram}

\label{Diagram Essential-Smooth-3-1-2}
\index{Diagram Essential-Smooth-3-1-2}

\centering
\setlength{\unitlength}{1mm}
{\renewcommand{\dashlinestretch}{30}
\begin{picture}(100,130)(0,0)
\put(50,80){\circle{80}}
\path(10,80)(90,80)

\put(100,100){{\xssc $X$}}
\put(100,60){{\xssc $\xbm (X)$}}

\path(50,100)(50,60)
\path(48.5,63)(50,60)(51.5,63)
\put(50,101){\circle*{0.3}}
\put(50,59){\circle*{0.3}}
\put(50,57){{\xssc $<x,i>$}}
\put(50,102){{\xssc $<y,j>$}}
\put(52,90){{\xssc $<\xba,k>$}}

\path(20,70)(48,70)
\path(46,71)(48,70)(46,69)
\put(19,70){\circle*{0.3}}
\put(13,67){{\xssc $<z',m'>$}}
\put(26,71){{\xssc $<\xbb',l'>$}}

\path(30,30)(30,68)
\path(29,66)(30,68)(31,66)
\put(30,29){\circle*{0.3}}
\put(31,53){{\xssc $<\xbg ',n'>$}}
\put(30,26.5){{\xssc $<u',p'>$}}

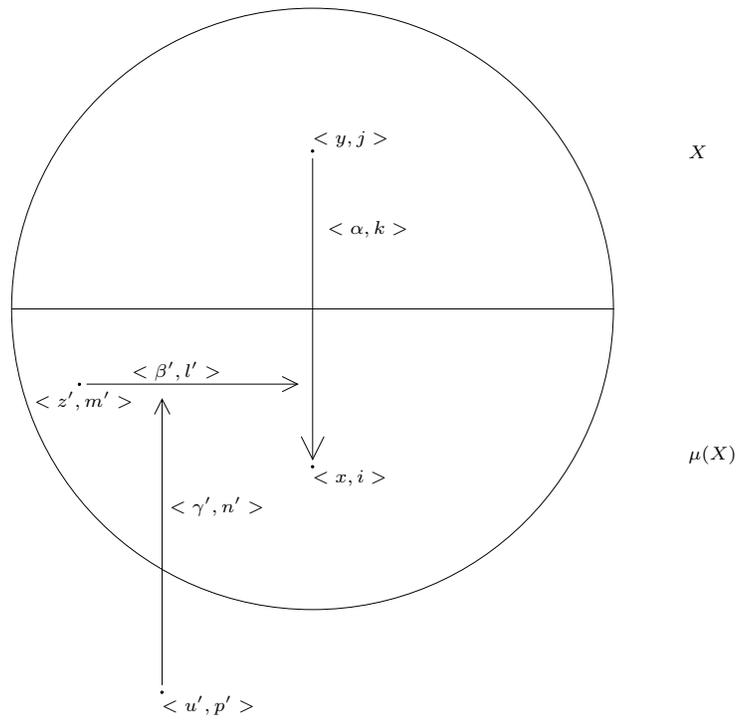
\put(30,10) {{\rm\bf Case 3-1-2}}

\end{picture}
}
\end{diagram}

\vspace{4mm}

\vspace{10mm}

\begin{diagram}

\label{Diagram Essential-Smooth-3-2}
\index{Diagram Essential-Smooth-3-2}

\centering
\setlength{\unitlength}{1mm}
{\renewcommand{\dashlinestretch}{30}
\begin{picture}(100,110)(0,0)
\put(50,60){\circle{80}}
\path(10,60)(90,60)

\put(100,80){{\xssc $X$}}
\put(100,40){{\xssc $\xbm (X)$}}

\path(50,80)(50,40)
\path(48.5,77)(50,80)(51.5,77)
\put(50,81){\circle*{0.3}}
\put(50,39){\circle*{0.3}}
\put(50,37){{\xssc $<y',j'>$}}
\put(50,82){{\xssc $<x,i>$}}
\put(52,70){{\xssc $<\xba',k'>$}}

\path(20,70)(48,70)
\path(46,71)(48,70)(46,69)
\put(19,70){\circle*{0.3}}
\put(13,67){{\xssc $<z',m'>$}}
\put(26,71){{\xssc $<\xbb',l'>$}}

\path(30,30)(30,68)
\path(29,66)(30,68)(31,66)
\put(30,29){\circle*{0.3}}
\put(31,53){{\xssc $<\xbg ',n'>$}}
\put(31,26.5){{\xssc $<u',p'>$}}

\put(30,0) {{\rm\bf Case 3-2}}

\end{picture}
}
\end{diagram}

\vspace{4mm}

\index{Fact X-Sub-X'}

\er

\bfa

$\hspace{0.01em}$


\label{Fact X-Sub-X'}

(1) If $X \xes X',$ then $X= \xbm (X' ),$

(2) $X \xes X',$ $X \xcc X'' \xcc X' $ $ \xch $ $X \xes X''.$ (This
corresponds to $( \xbm CUM).)$

(3) $X \xes X',$ $X \xcc Y',$ $Y \xes Y',$ $Y \xcc X' $ $ \xch $ $X=Y.$
(This corresponds to $( \xbm \xcc \xcd ).)$
\index{Fact X-Sub-X' Proof}

\efa

\paragraph{
Proof Fact $X-Sub-X' $
}

$\hspace{0.01em}$


\label{Section Proof Fact X-Sub-X'}

\subparagraph{
Proof
}

$\hspace{0.01em}$


(1) Trivial by Fact \ref{Fact Higher-Validity} (1).

(2)

We have to show

(a) $ \xcA x \xbe X'' -X$ $ \xcA <x,i>$ $ \xcE \xba:x' \xcp <x,i>( \xba $
is a valid $X \xch X'' $ arrow), and

(b) $ \xcA x \xbe X$ $ \xcE <x,i>$
$( \xcA \xba:x' \xcp <x,i>(x' \xbe X'' $ $ \xch $ $ \xcE \xbb:x'' \xcp
\xba.( \xbb $ is a valid $X \xch X'' $ arrow))).

Both follow from the corresponding condition for $X \xch X',$ the
restriction of the
universal quantifier, and Fact \ref{Fact Higher-Validity} (2).

(3)

Let $x \xbe X-$Y.

(a) By $x \xbe X \xes X',$ $ \xcE <x,i>$ s.t.
$( \xcA \xba:x' \xcp <x,i>(x' \xbe X' $ $ \xch $ $ \xcE \xbb:x'' \xcp
\xba.( \xbb $ is a valid $X \xch X' $ arrow))).

(b) By $x \xce Y \xes \xcE \xba_{1}:x' \xcp <x,i>$ $ \xba_{1}$ is a valid
$Y \xch Y' $ arrow, in
particular $x' \xbe Y \xcc X'.$ Moreover, $ \xbl ( \xba_{1})=1.$

So by (a) $ \xcE \xbb_{2}:x'' \xcp \xba.( \xbb_{2}$ is a valid $X \xch X'
$ arrow), in particular $x'' \xbe X \xcc Y',$
moreover $ \xbl ( \xbb_{2})=2.$

It follows by induction from the definition of valid $A \xch B$ arrows
that

$ \xcA n \xcE \xba_{2m+1},$ $ \xbl ( \xba_{2m+1})=2m+1,$ $ \xba_{2m+1}$ a
valid $Y \xch Y' $ arrow and

$ \xcA n \xcE \xbb_{2m+2},$ $ \xbl ( \xbb_{2m+2})=2m+2,$ $ \xbb_{2m+2}$ a
valid $X \xch X' $ arrow,

which is impossible, as $ \xdx $ is a structure of finite level.

$ \xcz $
\\[3ex]
\index{Definition Totally-Smooth}

\bd

$\hspace{0.01em}$


\label{Definition Totally-Smooth}

Let $ \xdx $ be a generalized preferential structure, $X \xcc \xdP ( \xdx
).$

$ \xdx $ is called totally smooth for $X$ iff

(1) $ \xcA \xba:x \xcp y \xbe \xdA ( \xdx )(O( \xba ) \xcv D( \xba ) \xcc
X$ $ \xch $ $ \xcE \xba ':x' \xcp y.x' \xbe \xbm (X))$

(2) if $ \xba $ is valid, then $ \xba ' $ must also be valid.

(y a point or an arrow).

If $ \xdy \xcc \xdP ( \xdx ),$ then $ \xdx $ is called $ \xdy -$totally
smooth iff for all $X \xbe \xdy $
$ \xdx $ is totally smooth for $X.$
\index{Example Totally-Smooth}

\ed

\be

$\hspace{0.01em}$


\label{Example Totally-Smooth}

$X:=\{ \xba:a \xcp b,$ $ \xba ':b \xcp c,$ $ \xba '':a \xcp c,$ $ \xbb
:b \xcp \xba ' \}$ is not totally smooth,

$X:=\{ \xba:a \xcp b,$ $ \xba ':b \xcp c,$ $ \xba '':a \xcp c,$ $ \xbb
:b \xcp \xba ',$ $ \xbb ':a \xcp \xba ' \}$ is totally smooth.
\index{Example Need-Smooth}

\ee

\be

$\hspace{0.01em}$


\label{Example Need-Smooth}

Consider $ \xba ':a \xcp b,$ $ \xba '':b \xcp c,$ $ \xba:a \xcp c,$ $
\xbb:a \xcp \xba.$

Then $ \xbm (\{a,b,c\})=\{a\},$ $ \xbm (\{a,c\})=\{a,c\}.$
Thus, $( \xbm CUM)$ does not hold in this structure.
Note that there is no valid arrow from $ \xbm (\{a,b,c\})$ to $c.$
\index{Definition Essentially-Smooth}

\ee

\bd

$\hspace{0.01em}$


\label{Definition Essentially-Smooth}

Let $ \xdx $ be a generalized preferential structure, $X \xcc \xdP ( \xdx
).$

$ \xdx $ is called essentially smooth for $X$ iff $ \xbm (X) \xes X.$

If $ \xdy \xcc \xdP ( \xdx ),$ then $ \xdx $ is called $ \xdy
-$essentially smooth iff for all $X \xbe \xdy $
$ \xbm (X) \xes X.$
\index{Example Total-vs-Essential}

\ed

\be

$\hspace{0.01em}$


\label{Example Total-vs-Essential}

It is easy to see that we can distinguish total and essential smoothness
in richer structures, as the following Example shows:

We add an accessibility relation $R,$ and consider only those models which
are accessible.

Let e.g. $a \xcp b \xcp <c,0>,$ $<c,1>,$ without transitivity. Thus, only
$c$ has two
copies. This structure is essentially smooth, but of course not totally
so.

Let now mRa, mRb, $mR<c,0>,$ $mR<c,1>,$ $m' Ra,$ $m' Rb,$ $m' R<c,0>.$

Thus, seen from $m,$ $ \xbm (\{a,b,c\})=\{a,c\},$ but seen from $m',$ $
\xbm (\{a,b,c\})=\{a\},$
but $ \xbm (\{a,c\})=\{a,c\},$ contradicting (CUM).

$ \xcz $
\\[3ex]

\vspace{7mm}


\vspace{7mm}

\subsubsection{
Results on not necessarily smooth structures
}


\label{Section Reac-GenPref-NonSmooth}
\index{Section Reac-GenPref-NonSmooth}
\index{Example Need-Pr}

\ee

\be

$\hspace{0.01em}$


\label{Example Need-Pr}

We show here $( \xbm \xcc )+( \xbm \xcc \xcd )+( \xbm CUM)+( \xbm RatM)+(
\xcs )$ $ \xcH $ $( \xbm PR).$

Let $U:=\{a,b,c\}.$ Let $ \xdy = \xdp (U).$ So $( \xcs )$ is trivially
satisfied.
Set $f(X):=X$ for all $X \xcc U$ except for $f(\{a,b\})=\{b\}.$ Obviously,
this cannot be represented by a preferential structure and $( \xbm PR)$ is
false
for $U$ and $\{a,b\}.$ But it satisfies $( \xbm \xcc ),$ $( \xbm CUM),$ $(
\xbm RatM).$ $( \xbm \xcc )$ is trivial.
$( \xbm CUM):$ Let $f(X) \xcc Y \xcc X.$ If $f(X)=X,$ we are done.
Consider $f(\{a,b\})=\{b\}.$ If
$\{b\} \xcc Y \xcc \{a,b\},$ then $f(Y)=\{b\},$ so we are done again. It
is shown in
Fact \ref{Fact Mu-Base}, (8) that $( \xbm \xcc \xcd )$ follows.
$( \xbm RatM):$ Suppose $X \xcc Y,$ $X \xcs f(Y) \xEd \xCQ,$ we have to
show $f(X) \xcc f(Y) \xcs X.$ If $f(Y)=Y,$ the result holds by $X \xcc Y,$
so it does if $X=Y.$
The only remaining case is $Y=\{a,b\},$ $X=\{b\},$ and the result holds
again.

$ \xcz $
\\[3ex]

\ee

The idea to solve the representation problem illustrated by Example \ref{Example
Need-Pr}
is to use the points $c$ and $d$ as bases for counterarguments against $
\xba:b \xcp a$ - as
is possible in IBRS. We do this now. We will obtain a representation
for logics weaker than $ \xCf P$ by generalized preferential structures.

\bd

$\hspace{0.01em}$


\label{Definition Legal-Sub}

NOETIG???

(1) Let $X$ be a diagram and $Y \xcc X.$
$Y$ is called a legal subdiagram of $X$ iff it can be obtained from $X$
and the set
of points of $Y$ inductively as follows:

(1.1) all points of $Y$ are in $Y$

(1.2) if $a,b$ (b a point or arrow) are in $Y,$ then so is any arrow from
a to $b,$
which is in $X.$

(These are all the objects in $Y).$

(2) If we work with copies of points, then we define:
$Y$ is called a legal subdiagram of $X$ iff it can be obtained from $X$
and the set
of points of $Y$ inductively as follows:

(2.1) if there is some $<y,i>$ in $Y \xcc X,$ then all $<y,i' > \xbe X$
are in $Y$ too,

(2.2) if $a,b$ (b a point or arrow) are in $Y,$ then so is any arrow from
a to $b,$
which is in $X.$

(These are all the objects in $Y).$

\ed

From now on, all subdiagrams will be legal.

We will now prove a representation theorem, but will make it more general
than
for preferential structures only. For this purpose, we will introduce some
definitions first.

\bd

$\hspace{0.01em}$


\label{Definition Eta-Rho-Structure}

Let $ \xbh, \xbr: \xdy \xcp \xdp (U).$

(1) If $ \xdx $ is a simple structure:

$ \xdx $ is called an attacking structure relative to $ \xbh $
representing $ \xbr $ iff

$ \xbr (X)$ $=$ $\{x \xbe \xbh (X):$ there is no valid $X-to- \xbh (X)$
arrow $ \xba:x' \xcp x\}$

for all $X \xbe \xdy.$

(2) If $ \xdx $ is a structure with copies:

$ \xdx $ is called an attacking structure relative to $ \xbh $
representing $ \xbr $ iff

$ \xbr (X)$ $=$ $\{x \xbe \xbh (X):$ there is $<x,i>$ and no valid $X-to-
\xbh (X)$ arrow
$ \xba:<x',i' > \xcp <x,i>\}$

for all $X \xbe \xdy.$

Obviously, in those cases $ \xbr (X) \xcc \xbh (X)$ for all $X \xbe \xdy
.$

Thus, $ \xdx $ is a preferential structure iff $ \xbh $ is the identity.

See Diagram \ref{Diagram Eta-Rho-1}

\vspace{10mm}

\begin{diagram}

\label{Diagram Eta-Rho-1}
\index{Diagram Eta-Rho-1}

\centering
\setlength{\unitlength}{1mm}
{\renewcommand{\dashlinestretch}{30}
\begin{picture}(150,100)(0,0)

\put(30,50){\circle{50}}
\put(100,50){\circle{40}}
\put(100,50){\circle{70}}

\put(30,50){\circle*{1}}
\put(70,50){\circle*{1}}

\path(31,50)(69,50)
\path(66.6,50.9)(69,50)(66.6,49.1)

\put(30,60){\xssc{$X$}}
\put(100,80){\xssc{$\xbh(X)$}}
\put(100,60){\xssc{$\xbr(X)$}}

\put(30,10) {{\rm\bf Attacking structure}}

\end{picture}
}

\end{diagram}

\vspace{4mm}

\ed

(Note that it does not seem very useful to generalize the notion of
smoothness
from preferential structures to general attacking structures, as, in the
general
case, the minimizing set $X$ and the result $ \xbr (X)$ may be disjoint.)

The following result is the first positive representation result of this
paper, and shows that we can obtain (almost) anything with level 2
structures.

\bp

$\hspace{0.01em}$


\label{Proposition Eta-Rho-Repres}

Let $ \xbh, \xbr: \xdy \xcp \xdp (U).$ Then there is an attacking level
2 structure relative
to $ \xbh $ representing $ \xbr $ iff

(1) $ \xbr (X) \xcc \xbh (X)$ for all $X \xbe \xdy,$

(2) $ \xbr ( \xCQ )= \xbh ( \xCQ )$ if $ \xCQ \xbe \xdy.$

\ep

(2) is, of course, void for preferential structures.

\subparagraph{
Proof:
}

$\hspace{0.01em}$


\label{Section Proof:}

(A) The construction

We make a two stage construction.

(A.1) Stage 1.

In stage one, consider (almost as usual)

$ \xdu:=< \xdx,\{ \xba_{i}:i \xbe I\}>$ where

$ \xdx $ $:=$ $\{<x,f>:$ $x \xbe U,$ $f \xbe \xbP \{X \xbe \xdy:x \xbe
\xbh (X)- \xbr (X)\}\},$

$ \xba:x' \xcp <x,f>$ $: \xcj $ $x' \xbe ran(f).$ Attention: $x' \xbe X,$
not $x' \xbe \xbr (X)!$

(A.2) Stage 2.

Let $ \xdx ' $ be the set of all $<x,f,X>$ s.t. $<x,f> \xbe \xdx $ and

(a) either $X$ is some dummy value, say $*$

or

(b) all of the following $(1)-(4)$ hold:

(1) $X \xbe \xdy,$

(2) $x \xbe \xbr (X),$

(3) there is $X' \xcc X,$ $x \xbe \xbh (X' )- \xbr (X' ),$ $X' \xbe \xdy
,$
(thus $ran(f) \xcs X \xEd \xCQ $ by definition),

(4) $ \xcA X'' \xbe \xdy.(X \xcc X'',$ $x \xbe \xbh (X'' )- \xbr (X'' )$
$ \xcp $ $(ran(f) \xcs X'' )-X \xEd \xCQ ).$

(Thus, $f$ chooses in (4) for $X'' $ also outside $X.$ If there is no such
$X'',$ (4) is
void, and only $(1)-(3)$ need to hold, i.e. we may take any $f$ with
$<x,f> \xbe \xdx.)$

See Diagram \ref{Diagram Condition-rho-eta}.

\vspace{10mm}

\begin{diagram}

\label{Diagram Condition-rho-eta}
\index{Diagram Condition-rho-eta}

\centering
\setlength{\unitlength}{0.00083333in}
{\renewcommand{\dashlinestretch}{30}
\begin{picture}(2727,2755)(0,-500)
\put(1304.562,-118.818){\arc{3584.794}{4.1827}{5.2252}}
\put(1356.377,-664.180){\arc{3978.201}{4.0248}{5.3839}}
\put(1370.659,430.666){\arc{1327.952}{3.9224}{5.4136}}
\put(1344,1343){\ellipse{2672}{2672}}
\put(1334,1293){\ellipse{1178}{1178}}
\put(1289,1338){\ellipse{1802}{1802}}
\put(2200,1950){\circle*{30}}
\put(2094,1800){{\xssc $f(X'')$}}
\put(2214,1318){{\xssc $\rho(X)$}}
\put(2074,2608){{\xssc $X''$}}
\put(1709,2228){{\xssc $X$}}
\put(2624,698) {{\xssc $\rho(X'')$}}
\put(1334,828) {{\xssc $\rho(X')$}}
\put(1404,1928){{\xssc $X'$}}
\put(1514,1433){\circle*{30}}
\put(1544,1403){{\xssc $x$}}
\put(100,2800){{\xssc For simplicity, $\xbh(X)=X$ here}}

\put(150,-400) {{\rm\bf The complicated case}}

\end{picture}
}
\end{diagram}

\vspace{4mm}

Note: If $(1)-(3)$ are satisfied for $x$ and $X,$ then we will find $f$
s.t. $<x,f> \xbe \xdx,$
and $<x,f,X>$ satisfies $(1)-(4):$ As $X \xcB X'' $ for $X'' $ as in (4),
we find $f$ which
chooses for such $X'' $ outside of $X.$

So for any $<x,f> \xbe \xdx,$ there is $<x,f,*>,$ and maybe also some
$<x,f,X>$ in $ \xdx '.$

Let again for any $x',$ $<x,f,X> \xbe \xdx ' $

$ \xba:x' \xcp <x,f,X>$ $: \xcj $ $x' \xbe ran(f)$

(A.3) Adding arrows.

Consider $x' $ and $<x,f,X>.$

If $X=*,$ or $x' \xce X,$ we do nothing, i.e. leave a simple arrow
$ \xba:x' \xcp <x,f,X>$ $ \xcj $ $x' \xbe ran(f).$

If $X \xbe \xdy,$ and $x' \xbe X,$ and $x' \xbe ran(f),$ we make $X$ many
copies of the attacking
arrow and have then: $< \xba,x'' >:x' \xcp <x,f,X>$ for all $x'' \xbe X.$

In addition, we add attacks on the $< \xba,x'' >:$ $< \xbb,x'' >:x''
\xcp < \xba,x'' >$ for all $x'' \xbe X.$

The full structure $ \xdz $ is thus:

$ \xdx ' $ is the set of elements.

If $x' \xbe ran(f),$ and $X=*$ or $x' \xce X$ then $ \xba:x' \xcp
<x,f,X>$

if $x' \xbe ran(f),$ and $X \xEd *$ and $x' \xbe X$ then

(a) $< \xba,x'' >:x' \xcp <x,f,X>$ for all $x'' \xbe X,$

(b) $< \xbb,x'' >:x'' \xcp < \xba,x'' >$ for all $x'' \xbe X.$

See Diagram \ref{Diagram Structure-rho-eta}.

\vspace{10mm}

\begin{diagram}

\label{Diagram Structure-rho-eta}
\index{Diagram Structure-rho-eta}

\centering
\setlength{\unitlength}{1mm}
{\renewcommand{\dashlinestretch}{30}
\begin{picture}(150,100)(0,0)

\put(30,50){\circle{50}}
\put(100,50){\circle{40}}
\put(100,50){\circle{70}}

\put(30,50){\circle*{1}}
\put(90,50){\circle*{1}}
\put(92,50){{\xssc $< x, f, X>$}}
\put(32,50) {{\xssc $x'$}}

\put(60,20){\arc{84.5}{-2.34}{-0.8}}
\put(70,62){{\xssc $< a, x_0''>$}}
\put(60,80){\arc{84.5}{0.8}{2.34}}
\put(70,38){{\xssc $< a, x_1''>$}}

\path(88.9,51.7)(90,50)(88,50.5)
\path(88,49.5)(90,50)(88.9,48.3)

\put(35,65){\circle*{1}}
\put(35,35){\circle*{1}}

\put(36,65) {{\xssc $x_0''$}}
\put(36,35){{\xssc $x''_1$}}

\path(35,65)(40,58)
\path(39.7,59.2)(40,58)(39.0,58.7)
\put(25,62){{\xssc $< \beta, x_1''>$}}

\path(35,35)(40,42)
\path(39,41.3)(40,42)(39.7,40.8)
\put(25,38) {{\xssc $< \beta, x_o''>$}}

\put(15,50){\xssc{$X$}}
\put(100,80){\xssc{$\xbh(X)$}}
\put(100,60){\xssc{$\xbr(X)$}}

\put(30,10) {{\rm\bf Attacking structure}}

\end{picture}
}

\end{diagram}

\vspace{4mm}

(B) Representation

We have to show that this structure represents $ \xbr $ relative to $ \xbh
.$

Let $y \xbe \xbh (Y),$ $Y \xbe \xdy.$

Case 1. $y \xbe \xbr (Y).$

We have to show that there is $<y,g,Y'' >$ s.t. there is no valid
$ \xba:y' \xcp <y,g,Y'' >,$ $y' \xbe Y.$ In Case 1.1 below, $Y'' $ will
be $*,$ in Case 1.2,
$Y'' $ will be $Y,$ $g$ will be chosen suitably.

Case 1.1. There is no $Y' \xcc Y,$ $y \xbe \xbh (Y' )- \xbr (Y' ),$ $Y'
\xbe \xdy.$

So for all $Y' $ with $y \xbe \xbh (Y' )- \xbr (Y' )$ $Y' -Y \xEd \xCQ.$
Let $g \xbe \xbP \{Y' -Y:y \xbe \xbh (Y' )- \xbr (Y' )\}.$
Then $ran(g) \xcs Y= \xCQ,$ and $<y,g,*>$ is not attacked from $Y.$
$(<y,g>$ was already not attacked in $ \xdx.)$

Case 1.2. There is $Y' \xcc Y,$ $y \xbe \xbh (Y' )- \xbr (Y' ),$ $Y' \xbe
\xdy.$

Let now $<y,g,Y> \xbe \xdx ',$ s.t. $g(Y'' ) \xce Y$ if $Y \xcc Y'',$ $y
\xbe \xbh (Y'' )- \xbr (Y'' ),$ $Y'' \xbe \xdy.$
As noted above, such $g$ and thus $<y,g,Y>$ exist. Fix $<y,g,Y>.$

Consider any $y' \xbe ran(g).$ If $y' \xce Y,$ $y' $ does not attack
$<y,g,Y>$ in $Y.$ Suppose $y' \xbe Y.$ We had made $Y$ many copies
$< \xba,y'' >,$ $y'' \xbe Y$ with $< \xba,y'' >:y' \xcp <y,g,Y>$ and had
added the level 2
arrows $< \xbb,y'' >:y'' \xcp < \xba,y'' >$ for $y'' \xbe Y.$
So all copies $< \xba,y'' >$ are destroyed in $Y.$ This was done for all
$y' \xbe Y,$ $y' \xbe ran(g),$ so $<y,g,Y>$ is now not (validly) attacked
in $Y$
any more.

Case 2. $y \xbe \xbh (Y)- \xbr (Y).$

Let $<y,g,Y' >$ (where $Y' $ can be $*)$ be any copy of $y,$ we have to
show that
there is $z \xbe Y,$ $ \xba:z \xcp <y,g,Y' >,$ or some
$< \xba,z' >:z \xcp <y,g,Y' >,$ $z' \xbe Y',$ which is not destroyed by
some
level 2 arrow $< \xbb,z' >:z' \xcp < \xba,z' >,$ $z' \xbe Y.$

As $y \xbe \xbh (Y)- \xbr (Y),$ $ran(g) \xcs Y \xEd \xCQ,$ so there is $z
\xbe ran(g) \xcs Y.$ Fix such $z.$
(We will modify the choice of $z$ only in Case 2.2.2 below.)

Case 2.1. $Y' =*.$

As $z \xbe ran(g),$ $ \xba:z \xcp <y,g,*>.$ (There were no level 2 arrows
introduced for
this copy.)

Case 2.2. $Y' \xEd *.$

So $<y,g,Y' >$ satisfies the conditions $(1)-(4)$ of (b) at the beginning
of the
proof.

If $z \xce Y',$ we are done, as $ \xba:z \xcp <y,g,Y' >,$ and there were
no level 2
arrows introduced in this case.
If $z \xbe Y',$ we had made $Y' $ many copies $< \xba,z' >,$ $< \xba,z'
>:z \xcp <y,g,Y' >,$ one for
each $z' \xbe Y'.$ Each $< \xba,z' >$ was destroyed by $< \xbb,z' >:z'
\xcp < \xba,z' >,$ $z' \xbe Y'.$

Case 2.2.1. $Y' \xcC Y.$

Let $z'' \xbe Y' -$Y, then $< \xba,z'' >:z \xcp <y,g,Y' >$ is destroyed
only by
$< \xbb,z'' >:z'' \xcp < \xba,z'' >$ in $Y',$ but not in $Y,$ as $z''
\xce Y,$ so $<y,g,Y' >$ is
attacked by $< \xba,z'' >:z \xcp <y,g,Y' >,$ valid in $Y.$

Case 2.2.2. $Y' \xcB Y$ $(Y=Y' $ is impossible, as $y \xbe \xbr (Y' ),$ $y
\xce \xbr (Y)).$

Then there was by definition (condition (b) (4)) some $z' \xbe (ran(g)
\xcs Y)-Y' $ and
$ \xba:z' \xcp <y,g,Y' >$ is valid, as $z' \xce Y'.$ (In this case,
there
are no copies of $ \xba $ and no level 2 arrows.)

$ \xcz $
\\[3ex]

\bco

$\hspace{0.01em}$


\label{Corollary Eta-Rho}

(1) We cannot distinguish general structures of level 2 from those of
higher
levels by their $ \xbr -$functions relative to $ \xbh.$

(2) Let $U$ be the universe, $ \xdy \xcc \xdp (U),$ $ \xbm: \xdy \xcp
\xdp (U).$
Then any $ \xbm $ satisfying $( \xbm \xcc )$ can be represented by a level
2 preferential
structure. (Choose $ \xbh =identity.)$

Again, we cannot distinguish general structures of level 2 from those of
higher
levels by their $ \xbm -$functions.

$ \xcz $
\\[3ex]

\eco

A remark on the function $ \xbh:$

We can also obtain the function $ \xbh $ via arrows. Of course, then we
need
positive arrows (not only negative arrows against negative arrows, as we
first
need to have something positive).

If $ \xbh $ is the identity, we can make a positive arrow from each point
to itself.
Otherwise, we can connect every point to every point by a positive arrow,
and then choose those we really want in $ \xbh $ by a choice function
obtained from
arrows just as we obtained $ \xbr $ from arrows.

\vspace{7mm}


\vspace{7mm}

\subsubsection{
Results on total smoothness
}


\label{Section Reac-GenPref-TotalSmooth}
\index{Section Reac-GenPref-TotalSmooth}

\bfa

$\hspace{0.01em}$


\label{Fact Val-Arrow}

Let $X,Y \xbe \xdy,$ $ \xdx $ a level $n$ structure.
Let $< \xba,k>:<x,i> \xcp <y,j>,$ where $<y,j>$ may itself be (a copy of)
an arrow.

(1) Let $n>1,$ $X \xcc Y,$ $< \xba,k> \xbe X$ a level $n-1$ arrow in $
\xdx \xex X.$ If $< \xba,k>$ is valid in
$ \xdx \xex Y,$ then it is valid in $ \xdx \xex X.$

(2) Let $ \xdx $ be totally smooth, $ \xbm (X) \xcc Y,$ $ \xbm (Y) \xcc
X,$ $< \xba,k> \xbe \xdx \xex X \xcs Y,$ then $< \xba,k>$
is valid in $ \xdx \xex X$ iff it is valid in $ \xdx \xex Y.$

Note that we will also sometimes write $X$ for $ \xdx \xex X,$ when the
context is clear.

\efa

\subparagraph{
Proof:
}

$\hspace{0.01em}$


\label{Section Proof:}

(1)
If $< \xba,k>$ is not valid in $ \xdx \xex X,$ then there must be a level
$n$ arrow
$< \xbb,r>:<z,s> \xcp < \xba,k>$ in $ \xdx \xex X \xcc \xdx \xex Y.$ $<
\xbb,r>$ must be valid in $ \xdx \xex X$ and $ \xdx \xex Y,$
as there are no level $n+1$ arrows. So $< \xba,k>$ is not valid in $ \xdx
\xex Y,$ $contradiction.$

(2)
By downward induction.
Case $n:$ $< \xba,k> \xbe \xdx \xex X \xcs Y,$ so it is valid in both as
there are no level $n+1$
arrows.
Case $m \xcp m-1:$ Let $< \xba,k> \xbe \xdx \xex X \xcs Y$ be a level
$m-1$ arrow valid in $ \xdx \xex X,$ but not
in $ \xdx \xex Y.$ So there must be a level $m$ arrow $< \xbb,r>:<z,s>
\xcp < \xba,k>$ valid in $ \xdx \xex Y.$
By total smoothness, we may assume $z \xbe \xbm (Y) \xcc X,$ so $< \xbb
,r> \xbe \xdx \xex X$ valid by
induction hypothesis. So $< \xba,k>$ is not valid in $ \xdx \xex X,$
$contradiction.$

$ \xcz $
\\[3ex]

\bco

$\hspace{0.01em}$


\label{Corollary Total-Mu}

Let $X,Y \xbe \xdy,$ $ \xdx $ a totally smooth level $n$ structure, $
\xbm (X) \xcc Y,$ $ \xbm (Y) \xcc X.$
Then $ \xbm (X)= \xbm (Y).$

\eco

\subparagraph{
Proof:
}

$\hspace{0.01em}$


\label{Section Proof:}

Let $x \xbe \xbm (X)- \xbm (Y).$ Then by $ \xbm (X) \xcc Y,$ $x \xbe Y,$
so there must be for all $<x,i> \xbe \xdx $
an arrow $< \xba,k>:<y,j> \xcp <x,i>$ valid in $ \xdx \xex Y,$ wlog. $y
\xbe \xbm (Y) \xcc X$ by total
smoothness. So by Fact \ref{Fact Val-Arrow}, (2),
$< \xba,k>$ is valid in $ \xdx \xex X.$ This
holds for all $<x,i>,$ so $x \xce \xbm (X),$ $contradiction.$
$ \xcz $
\\[3ex]

\bfa

$\hspace{0.01em}$


\label{Fact Level-bigger-2}

There are situations satisfying $( \xbm \xcc )+( \xbm CUM)+( \xcs )$ which
cannot be represented
by level 2 totally smooth preferential structures.

\efa

The proof is given in the following example.

\be

$\hspace{0.01em}$


\label{Example Level-bigger-2}

Let $Y:=\{x,y,y' \},$ $X:=\{x,y\},$ $X':=\{x,y' \}.$ Let $ \xdy:= \xdp
(Y).$
Let $ \xbm (Y):=\{y,y' \},$ $ \xbm (X):= \xbm (X' ):=\{x\},$ and $ \xbm
(Z):=Z$ for all other sets.

Obviously, this satisfies $( \xcs ),$ $( \xbm \xcc ),$ and $( \xbm CUM).$

Suppose $ \xdx $ is a totally smooth level 2 structure representing $ \xbm
.$

So $ \xbm (X)= \xbm (X' ) \xcc Y- \xbm (Y),$ $ \xbm (Y) \xcc X \xcv X'.$
Let $<x,i>$ be minimal in $ \xdx \xex X.$
As $<x,i>$ cannot be minimal in $ \xdx \xex Y,$ there must be $ \xba
:<z,j> \xcp <x,i>,$ valid in
$ \xdx \xex Y.$

Case 1: $z \xbe X'.$

So $ \xba \xbe \xdx \xex X'.$
If $ \xba $ is valid in $ \xdx \xex X',$ there must be $ \xba ':<x',i'
> \xcp <x,i>,$ $x' \xbe \xbm (X' ),$ valid in
$ \xdx \xex X',$ and thus in $ \xdx \xex X,$ by $ \xbm (X)= \xbm (X' )$
and
Fact \ref{Fact Val-Arrow} (2). This is
impossible, so there must be $ \xbb:<x',i' > \xcp \xba,$ $x' \xbe \xbm
(X' ),$ valid in $ \xdx \xex X'.$
As $ \xbb $ is in $ \xdx \xex Y$ and $ \xdx $ a level $ \xck 2$ structure,
$ \xbb $ is valid in $ \xdx \xex Y,$ so
$ \xba $ is not valid in $ \xdx \xex Y,$ $contradiction.$

Case 2: $z \xbe X.$

$ \xba $ cannot be valid in $ \xdx \xex X,$ so there must be $ \xbb:<x'
,i' > \xcp \xba,$ $x' \xbe \xbm (X),$
valid in $ \xdx \xex X.$ Again, as $ \xbb $ is in $ \xdx \xex Y$ and $
\xdx $ a level $ \xck 2$ structure, $ \xbb $ is
valid in $ \xdx \xex Y,$ so $ \xba $ is not valid in $ \xdx \xex Y,$
$contradiction.$

$ \xcz $
\\[3ex]

\ee

It is unknown to the authors whether an analogon is true for essential
smoothness, i.e. whether there are examples of such $ \xbm $ function
which need
at least level 3 essentially smooth structures for representation.
Proposition \ref{Proposition Level-3-Repr} below shows
that such structures suffice, but we do not
know whether level 3 is necessary.

\bfa

$\hspace{0.01em}$


\label{Fact Level-3-Solution}

Above Example \ref{Example Level-bigger-2}
can be solved by a totally smooth level 3 structure:

Let $ \xba_{1}:x \xcp y,$ $ \xba_{2}:x \xcp y',$ $ \xba_{3}:y \xcp x,$
$ \xbb_{1}:y \xcp \xba_{2},$ $ \xbb_{2}:y' \xcp \xba_{1},$ $ \xbb_{3}:y
\xcp \xba_{3},$ $ \xbb_{4}:x \xcp \xba_{3},$
$ \xbg_{1}:y' \xcp \xbb_{3},$ $ \xbg_{2}:y' \xcp \xbb_{4}.$

See Diagram \ref{Diagram Smooth-Level-3}.

\vspace{10mm}

\begin{diagram}

\label{Diagram Smooth-Level-3}
\index{Diagram Smooth-Level-3}

\centering
\setlength{\unitlength}{0.00083333in}
{\renewcommand{\dashlinestretch}{30}
\begin{picture}(3584,3681)(0,-500)
\put(1980.858,1684.979){\arc{2896.137}{3.0863}{5.1276}}
\path(503.037,1724.492)(535.000,1605.000)(563.029,1725.476)
\put(613.871,2084.678){\arc{1045.415}{1.8998}{5.0662}}
\path(671.312,2573.903)(795.000,2575.000)(685.347,2632.239)
\put(1755.902,1748.510){\arc{3480.502}{1.0969}{3.3803}}
\path(69.787,2036.400)(65.000,2160.000)(11.059,2048.688)
\put(1915.859,2762.471){\arc{1397.908}{3.4506}{5.9440}}
\path(1268.561,3097.293)(1250.000,2975.000)(1323.927,3074.171)
\put(1892.516,1773.709){\arc{3365.589}{4.8901}{7.4299}}
\path(2313.635,3433.788)(2190.000,3430.000)(2300.872,3375.161)
\path(605,1595)(2535,2930)
\path(2453.376,2837.062)(2535.000,2930.000)(2419.243,2886.407)
\path(600,1540)(2545,245)
\put(550,1560){\circle*{30}}
\path(2428.488,286.533)(2545.000,245.000)(2461.741,336.476)
\path(2560,2870)(1820,765)
\path(1831.496,888.158)(1820.000,765.000)(1888.100,868.259)
\path(2530,270)(1745,2345)
\put(2570,220){\circle*{30}}
\put(2600,2950){\circle*{30}}
\path(1815.520,2243.378)(1745.000,2345.000)(1759.401,2222.148)
\put(3315,2745){{\xssc $\gamma_1$}}
\put(2570,2845){{\xssc $y$}}
\put(2575,60)  {{\xssc $y'$}}
\put(1960,950) {{\xssc $\beta_1$}}
\put(1915,1985){{\xssc $\beta_2$}}
\put(2030,2445){{\xssc $\alpha_1$}}
\put(520,1435) {{\xssc $x$}}
\put(1960,3510){{\xssc $\beta_3$}}
\put(750,185)  {{\xssc $\gamma_2$}}
\put(170,2530) {{\xssc $\beta_4$}}
\put(630,2110) {{\xssc $\alpha_3$}}
\put(1920,500) {{\xssc $\alpha_2$}}

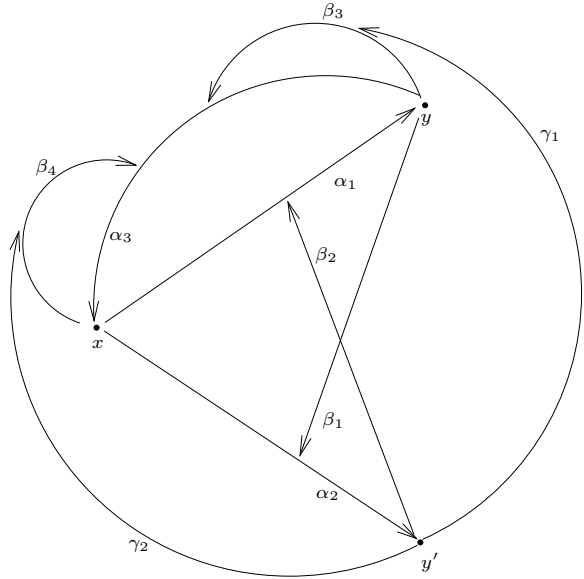
\put(100,-400) {{\rm\bf Solution by smooth level 3 structure}}

\end{picture}
}
\end{diagram}

\vspace{4mm}

\efa

The legal subdiagram generated by $X$ contains $ \xba_{1},$ $ \xba_{3},$ $
\xbb_{3},$ $ \xbb_{4}.$
$ \xba_{1},$ $ \xbb_{3},$ $ \xbb_{4}$ are valid, so $ \xbm (X)=\{x\}.$

The legal subdiagram generated by $X' $ contains $ \xba_{2}.$
$ \xba_{2}$ is valid, so $ \xbm (X' )=\{x\}.$

In the full diagram, $ \xba_{3},$ $ \xbb_{1},$ $ \xbb_{2},$ $ \xbg_{1},$ $
\xbg_{2}$ are valid, so $ \xbm (Y)=\{y,y' \}.$

$ \xcz $
\\[3ex]

\br

$\hspace{0.01em}$


\label{Remark Need-Mucd}

Example \ref{Example Mu-Cum-Cd} together with Corollary \ref{Corollary Total-Mu}
show
that $( \xbm \xcc )$ and $( \xbm CUM)$ without $( \xcs )$ do not guarantee
representability by a
level $n$ totally smooth structure.
\index{Example Mu-Cum-Cd}

\er

\be

$\hspace{0.01em}$


\label{Example Mu-Cum-Cd}

We show here $( \xbm \xcc )+( \xbm CUM)$ $ \xcH $ $( \xbm \xcc \xcd ).$

Consider $X:=\{a,b,c\},$ $Y:=\{a,b,d\},$ $f(X):=\{a\},$ $f(Y):=\{a,b\},$ $
\xdy:=\{X,Y\}.$
(If $f(\{a,b\})$ were defined, we would have $f(X)=f(\{a,b\})=f(Y),$
$contradiction.)$

Obviously, $( \xbm \xcc )$ and $( \xbm CUM)$ hold, but not $( \xbm \xcc
\xcd ).$

$ \xcz $
\\[3ex]

\vspace{7mm}


\vspace{7mm}

\subsubsection{
Results on essential smoothness
}


\label{Section Reac-GenPref-EssSmooth}
\index{Section Reac-GenPref-EssSmooth}

\ee

\bd

$\hspace{0.01em}$


\label{Definition ODPi}

Let $ \xbm: \xdy \xcp \xdp (U)$ and $ \xdx $ be given, let $ \xba:<y,j>
\xcp <x,i> \xbe \xdx.$

Define

$ \xdO ( \xba )$ $:=$ $\{Y \xbe \xdy:x \xbe Y- \xbm (Y),y \xbe \xbm
(Y)\},$

$ \xdD ( \xba )$ $:=$ $\{X \xbe \xdy:x \xbe \xbm (X),y \xbe X\},$

$ \xbP ( \xdO, \xba )$ $:=$ $ \xbP \{ \xbm (Y):Y \xbe \xdO ( \xba )\},$

$ \xbP ( \xdD, \xba )$ $:=$ $ \xbP \{ \xbm (X):X \xbe \xdD ( \xba )\}.$

\ed

\bl

$\hspace{0.01em}$


\label{Lemma Level-3-Constr}

Let $U$ be the universe, $ \xbm: \xdy \xcp \xdp (U).$
Let $ \xbm $ satisfy $( \xbm \xcc )+( \xbm \xcc \xcd ).$

Let $ \xdx $ be a level 1 preferential structure, $ \xba:<y,j> \xcp
<x,i>,$ $ \xdO ( \xba ) \xEd \xCQ,$
$ \xdD ( \xba ) \xEd \xCQ.$

We can modify $ \xdx $ to a level 3 structure $ \xdx ' $ by introducing
level 2 and level 3
arrows s.t. no copy of $ \xba $ is valid in any $X \xbe \xdD ( \xba ),$
and in every $Y \xbe \xdO ( \xba )$ at
least one copy of $ \xba $ is valid. (More precisely, we should write $
\xdx ' \xex X$ etc.)

Thus, in $ \xdx ',$

(1) $<x,i>$ will not be minimal in any $Y \xbe \xdO ( \xba ),$

(2) if $ \xba $ is the only arrow minimizing $<x,i>$ in $X \xbe \xdD (
\xba ),$ $<x,i>$ will now be
minimal in $X.$

The construction is made independently for all such arrows $ \xba \xbe
\xdx.$

(This is probably the main technical result of the paper.)

\el

\subparagraph{
Proof:
}

$\hspace{0.01em}$


\label{Section Proof:}

(1) The construction

Make $ \xbP ( \xdD, \xba )$ many copies of $ \xba:$ $\{< \xba,f>:f \xbe
\xbP ( \xdD, \xba )\},$ all
$< \xba,f>:<y,j> \xcp <x,i>.$ Note that $< \xba,f> \xbe X$ for all $X
\xbe \xdD ( \xba )$ and $< \xba,f> \xbe Y$ for all
$Y \xbe \xdO ( \xba ).$

Add to the structure $< \xbb,f,X_{r},g>:<f(X_{r}),i_{r}> \xcp < \xba
,f>,$ for any $X_{r} \xbe \xdD ( \xba ),$ and
$g \xbe \xbP ( \xdO, \xba )$ (and some or all $i_{r}$ - this does not
matter).

For all $Y_{s} \xbe \xdO ( \xba ):$

if $ \xbm (Y_{s}) \xcC X_{r}$ and $f(X_{r}) \xbe Y_{s},$ then add to the
structure
$< \xbg,f,X_{r},g,Y_{s}>:<g(Y_{s}),j_{s}> \xcp < \xbb,f,X_{r},g>$ (again
for all or some $j_{s}),$

if $ \xbm (Y_{s}) \xcc X_{r}$ or $f(X_{r}) \xce Y_{s},$ $< \xbg
,f,X_{r},g,Y_{s}>$ is not added.

See Diagram \ref{Diagram Essential-Smooth-Repr}.

\vspace{30mm}

\begin{diagram}

\label{Diagram Essential-Smooth-Repr}
\index{Diagram Essential-Smooth-Repr}

\centering
\setlength{\unitlength}{1mm}
{\renewcommand{\dashlinestretch}{30}
\begin{picture}(110,150)(0,0)
\put(50,80){\ellipse{120}{80}}
\put(50,100){\ellipse{120}{80}}
\path(20,134.3)(20,65.5)
\path(80,114.3)(80,45.5)
\put(50,142){\xssc{$X \xbe \xdD(\xba)$}}
\put(50,36){\xssc{$Y \xbe \xdO(\xba)$}}
\put(10,138){\xssc{$\xbm(X)$}}
\put(90,44){\xssc{$\xbm(Y)$}}

\path(10,82)(90,82)
\put(9.2,82){\circle*{0.3}}
\put(90.8,82){\circle*{0.3}}
\path(12.8,81)(10,82)(12.8,83)
\put(5,78){{\xssc $<x,i>$}}
\put(85,78){{\xssc $<y,j>$}}
\put(55,78){{\xssc $<\xba,f>$}}

\path(8,100)(35,82)
\path(32.5,82.6)(35,82)(33.5,84.2)
\put(7.2,100.5){\circle*{0.3}}
\put(7.7,101){{\xssc $<f(X_r),i_r>$}}
\put(4,90){{\xssc $<\xbb,f,X_r,g>$}}

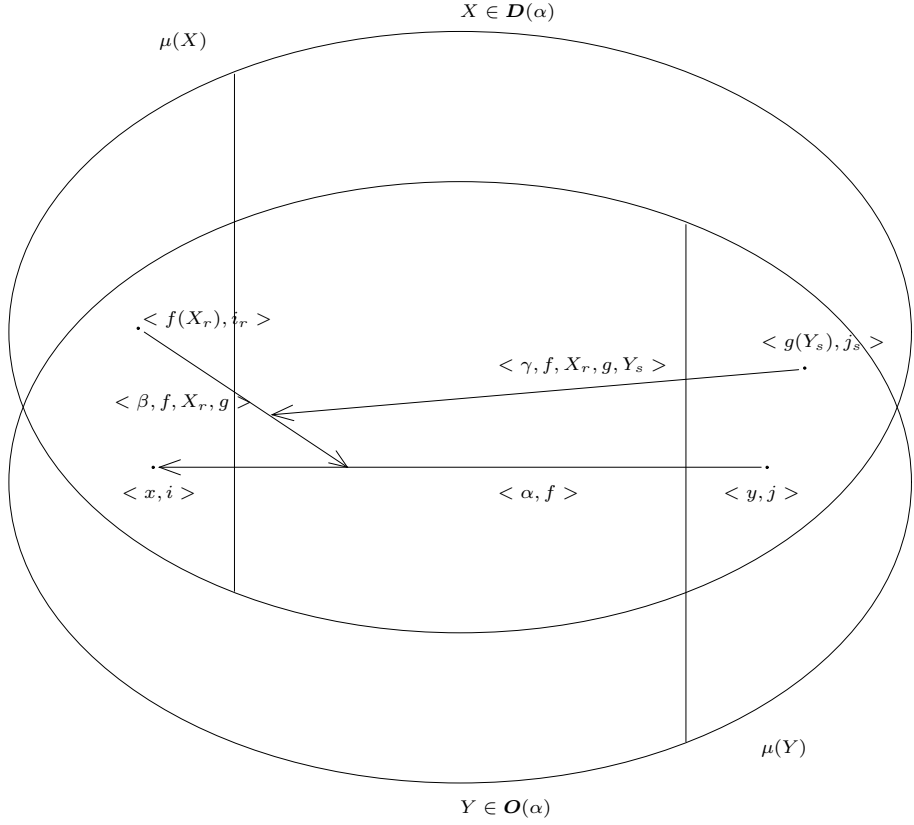
\path(25,89)(95,95)
\path(27.9,88.2)(25,89)(27.7,90.26)
\put(95.8,95.2){\circle*{0.3}}
\put(90,98){{\xssc $<g(Y_s),j_s>$}}
\put(55,95){{\xssc $<\xbg,f,X_r,g,Y_s>$}}

\put(10,15) {{\rm\bf The construction}}

\end{picture}
}

\end{diagram}

\vspace{4mm}

(2)

Let $X_{r} \xbe \xdD ( \xba ).$ We have to show that no $< \xba,f>$ is
valid in $X_{r}.$ Fix $f.$

$< \xba,f>$ is in $X_{r},$ so we have to show that for at least one $g
\xbe \xbP ( \xdO, \xba )$
$< \xbb,f,X_{r},g>$ is valid in $X_{r},$ i.e. that for this $g,$ no
$< \xbg,f,X_{r},g,Y_{s}>:<g(Y_{s}),j_{s}> \xcp < \xbb,f,X_{r},g>,$
$Y_{s} \xbe \xdO ( \xba )$ attacks $< \xbb,f,X_{r},g>$ in $X_{r}.$

We define $g.$ Take $Y_{s} \xbe \xdO ( \xba ).$

Case 1: $ \xbm (Y_{s}) \xcc X_{r}$ or $f(X_{r}) \xce Y_{s}:$ choose
arbitrary $g(Y_{s}) \xbe \xbm (Y_{s}).$

Case 2: $ \xbm (Y_{s}) \xcC X_{r}$ and $f(X_{r}) \xbe Y_{s}:$ Choose
$g(Y_{s}) \xbe \xbm (Y_{s})-X_{r}.$

In Case 1, $< \xbg,f,X_{r},g,Y_{s}>$ does not exist, so it cannot attack
$< \xbb,f,X_{r},g>.$

In Case 2, $< \xbg,f,X_{r},g,Y_{s}>:<g(Y_{s}),j_{s}> \xcp < \xbb
,f,X_{r},g>$ is not in $X_{r},$ as $g(Y_{s}) \xce X_{r}.$

Thus, no $< \xbg,f,X_{r},g,Y_{s}>:<g(Y_{s}),j_{s}> \xcp < \xbb
,f,X_{r},g>,$ $Y_{s} \xbe \xdO ( \xba )$ attacks $< \xbb,f,X_{r},g>$
in $X_{r}.$

So $ \xcA < \xba,f>:<y,j> \xcp <x,i>$

$ \xDC \xDC y \xbe X_{r}$ $ \xcp $ $ \xcE < \xbb
,f,X_{r},g>:<f(X_{r}),i_{r}> \xcp < \xba,f>$

$ \xDC \xDC \xDC \xDC (f(X_{r}) \xbe \xbm (X_{r})$ $ \xcu $ $ \xCN \xcE <
\xbg,f,X_{r},g,Y_{s}>:<g(Y_{s}),j_{s}> \xcp < \xbb,f,X_{r},g>.g(Y_{s})
\xbe X_{r}).$

But $< \xbb,f,X_{r},g>$ was constructed only for $< \xba,f>,$ so was $<
\xbg,f,X_{r},g,Y_{s}>,$ and
there was no other $< \xbg,i>$ attacking $< \xbb,f,X_{r},g>,$ so we are
done.

(3)

Let $Y_{s} \xbe \xdO ( \xba ).$ We have to show that at least one $< \xba
,f>$ is valid in $Y_{s}.$

We define $f \xbe \xbP ( \xdD, \xba ).$ Take $X_{r}.$

If $ \xbm (X_{r}) \xcC Y_{s},$ choose $f(X_{r}) \xbe \xbm (X_{r})-Y_{s}.$
If $ \xbm (X_{r}) \xcc Y_{s},$ choose arbitrary $f(X_{r}) \xbe \xbm
(X_{r}).$

All attacks on $<x,f>$ have the form $< \xbb,f,X_{r},g>:<f(X_{r}),i_{r}>
\xcp < \xba,f>,$ $X_{r} \xbe \xdD ( \xba ),$
$g \xbe \xbP ( \xdO, \xba ).$ We have to show that they are either not in
$Y_{s},$ or that they
are themselves attacked in $Y_{s}.$

Case 1: $ \xbm (X_{r}) \xcC Y_{s}.$ Then $f(X_{r}) \xce Y_{s},$ so $< \xbb
,f,X_{r},g>:<f(X_{r}),i_{r}> \xcp < \xba,f>$ is not
in $Y_{s}$ (for no $g).$

Case 2: $ \xbm (X_{r}) \xcc Y_{s}.$ Then $ \xbm (Y_{s}) \xcC X_{r}$ by $(
\xbm \xcc \xcd )$ and $f(X_{r}) \xbe Y_{s},$ so
$< \xbb,f,X_{r},g>:<f(X_{r}),i_{r}> \xcp < \xba,f>$ is in $Y_{s}$ (for
all $g).$ Take any $g \xbe \xbP ( \xdO, \xba ).$ As
$ \xbm (Y_{s}) \xcC X_{r}$ and $f(X_{r}) \xbe Y_{s},$ $< \xbg
,f,X_{r},g,Y_{s}>:<g(Y_{s}),j_{s}> \xcp < \xbb,f,X_{r},g>$ is defined,
and
$g(Y_{s}) \xbe \xbm (Y_{s}),$ so it is in $Y_{s}$ (for all $g).$ Thus, $<
\xbb,f,X_{r},g>$ is attacked in $Y_{s}.$

Thus, for this $f,$ all $< \xbb,f,X_{r},g>$ are either not in $Y_{s},$ or
attacked in $Y_{s},$
thus for this $f,$ $< \xba,f>$ is valid in $Y_{s}.$

So for this $<x,i>$

$ \xcE < \xba,f>:<y,j> \xcp <x,i>.y \xbe \xbm (Y_{s})$ $ \xcu $

$ \xDC \xDC $ (a) $ \xCN \xcE < \xbb,f,X_{r},g>:<f(X_{r}),i> \xcp < \xba
,f>.f(X_{r}) \xbe Y_{s}$

$ \xDC \xDC \xDC $ or

$ \xDC \xDC $ (b) $ \xcA < \xbb,f,X_{r},g>:<f(X_{r}),i> \xcp < \xba,f>$

$ \xDC \xDC \xDC \xDC (f(X_{r}) \xbe Y_{s}$ $ \xcp $

$ \xDC \xDC \xDC \xDC \xcE < \xbg,f,X_{r},g,Y_{s}>:<g(Y_{s}),j_{s}> \xcp
< \xbb,f,X_{r},g>.g(Y_{s}) \xbe \xbm (Y_{s})).$

As we made copies of $ \xba $ only, introduced only $ \xbb ' $s attacking
the $ \xba -$copies, and
$ \xbg ' $s attacking the $ \xbb ' $s, the construction is independent for
different $ \xba ' $s.

$ \xcz $
\\[3ex]

\bp

$\hspace{0.01em}$


\label{Proposition Level-3-Repr}

Let $U$ be the universe, $ \xbm: \xdy \xcp \xdp (U).$

Then any $ \xbm $ satisfying $( \xbm \xcc ),$ $( \xcs ),$ $( \xbm CUM)$
(or, alternatively, $( \xbm \xcc )$ and
$( \xbm \xcc \xcd ))$ can be represented by a level 3 essentially smooth
structure.

\ep

\subparagraph{
Proof:
}

$\hspace{0.01em}$


\label{Section Proof:}

In stage one, consider as usual $ \xdu:=< \xdx,\{ \xba_{i}:i \xbe I\}>$
where
$ \xdx $ $:=$ $\{<x,f>:$ $x \xbe U,$ $f \xbe \xbP \{ \xbm (X):$ $X \xbe
\xdy,$ $x \xbe X- \xbm (X)\}\},$ and set
$ \xba:<x',f' > \xcp <x,f>$ $: \xcj $ $x' \xbe ran(f).$

For stage two:

Any level 1 arrow $ \xba:<y,j> \xcp <x,i>$ was introduced in stage one by
some $Y \xbe \xdy $
s.t. $y \xbe \xbm (Y),$ $x \xbe Y- \xbm (Y).$ Do the construction
of Lemma \ref{Lemma Level-3-Constr} for all
level 1 arrows of $ \xdx $ in parallel or successively.

We have to show that the resulting structure represents $ \xbm $ and is
essentially smooth. (Level 3 is obvious.)

(1) Representation

Suppose $x \xbe Y- \xbm (Y).$ Then there was in stage 1 for all $<x,i>$
some $ \xba:<y,j> \xcp <x,i>,$
$y \xbe \xbm (Y).$ We examine the $y.$

If there is no $X$ s.t. $x \xbe \xbm (X),$ $y \xbe X,$ then there were no
$ \xbb ' $s and $ \xbg ' $s
introduced for this $ \xba:<y,j> \xcp <x,i>,$ so $ \xba $ is valid.

If there is $X$ s.t. $x \xbe \xbm (X),$ $y \xbe X,$ consider $ \xba:<y,j>
\xcp <x,i>.$ So $X \xbe \xdD ( \xba ),$
$Y \xbe \xdO ( \xba ),$ so we did the construction
of Lemma \ref{Lemma Level-3-Constr}, and by its
result, $<x,i>$ is not minimal in $Y.$

Thus, in both cases, $<x,i>$ is successfully attacked in $Y,$ and no
$<x,i>$
is a minimal element in $Y.$

Suppose $x \xbe \xbm (X)$ (we change notation to conform
to Lemma \ref{Lemma Level-3-Constr}).
Fix $<x,i>.$

If there is no $ \xba:<y,j> \xcp <x,i>,$ $y \xbe X,$ then $<x,i>$ is
minimal in $X,$ and we are
done.

If there is $ \xba $ or $< \xba,k>:<y,j> \xcp <x,i>,$ $y \xbe X,$ then $
\xba $ originated from stage one
through some $Y$ s.t. $x \xbe Y- \xbm (Y),$ and $y \xbe \xbm (Y).$ (Note
that stage 2 of the
construction did not
introduce any new level 1 arrows - only copies of existing level 1
arrows.) So
$X \xbe \xdD ( \xba ),$ $Y \xbe \xdO ( \xba ),$ so we did the construction
of Lemma \ref{Lemma Level-3-Constr}, and
by its result, $<x,i>$ is minimal in $X,$ and we are done again.

In both cases, all $<x,i>$ are minimal elements in $X.$

(2) Essential smoothness. We have to show the conditions of
Definition \ref{Definition X-Sub-X'}. We will, however, work with the
reformulation given in Remark \ref{Remark X-Sub-X'}.

Case (1), $x \xbe \xbm (X).$

Case (1.1), there is $<x,i>$ with no $< \xba,f>:<y,j> \xcp <x,i>,$ $y
\xbe X.$
There is nothing to show.

Case (1.2), for all $<x,i>$ there is $< \xba,f>:<y,j> \xcp <x,i>,$ $y
\xbe X.$

$ \xba $ was introduced in stage 1 by some $Y$ s.t.
$x \xbe Y- \xbm (Y),$ $y \xbe X \xcs \xbm (Y),$ so $X \xbe \xdD ( \xba ),$
$Y \xbe \xdO ( \xba ).$
In the proof of Lemma \ref{Lemma Level-3-Constr}, at the end of (2),
it was shown that

$ \xDC \xDC \xcE < \xbb,f,X_{r},g>:<f(X_{r}),i_{r}> \xcp < \xba,f>$

$ \xDC \xDC \xDC \xDC (f(X_{r}) \xbe \xbm (X_{r})$ $ \xcu $

$ \xDC \xDC \xDC \xDC \xCN \xcE < \xbg,f,X_{r},g,Y_{s}>:<g(Y_{s}),j_{s}>
\xcp < \xbb,f,X_{r},g>.g(Y_{s}) \xbe X_{r}).$

By $f(X_{r}) \xbe \xbm (X_{r}),$ condition (1) of Remark \ref{Remark X-Sub-X'}
is true.

Case (2), $x \xce \xbm (Y).$ Fix $<x,i>.$ (We change notation back to
$Y.)$

In stage 1, we constructed $ \xba:<y,j> \xcp <x,i>,$ $y \xbe \xbm (Y),$
so $Y \xbe \xdO ( \xba ).$

If $ \xdD ( \xba )= \xCQ,$ then there is no attack on $ \xba,$ and the
condition (2) of
Remark \ref{Remark X-Sub-X'} is trivially true.

If $ \xdD ( \xba ) \xEd \xCQ,$ we did the construction
of Lemma \ref{Lemma Level-3-Constr}, so

$ \xcE < \xba,f>:<y,j> \xcp <x,i>.y \xbe \xbm (Y_{s})$ $ \xcu $

$ \xDC \xDC $ (a) $ \xCN \xcE < \xbb,f,X_{r},g>:<f(X_{r}),i> \xcp < \xba
,f>.f(X_{r}) \xbe Y_{s}$

$ \xDC \xDC \xDC $ or

$ \xDC \xDC $ (b) $ \xcA < \xbb,f,X_{r},g>:<f(X_{r}),i> \xcp < \xba,f>$

$ \xDC \xDC \xDC \xDC (f(X_{r}) \xbe Y_{s}$ $ \xcp $

$ \xDC \xDC \xDC \xDC \xcE < \xbg,f,X_{r},g,Y_{s}>:<g(Y_{s}),j_{s}> \xcp
< \xbb,f,X_{r},g>.g(Y_{s}) \xbe \xbm (Y_{s}).$

As the only attacks on $< \xba,f>$ had the form $< \xbb,f,X_{r},g>,$ and
$g(Y_{s}) \xbe \xbm (Y_{s}),$
condition (2) of Remark \ref{Remark X-Sub-X'} is satisfied.

$ \xcz $
\\[3ex]

As said after Example \ref{Example Level-bigger-2},
we do not know if level 3 is necessary for
representation. We also do not know whether the same can be achieved with
level 3, or higher, totally smooth structures.

\vspace{7mm}


\vspace{7mm}


\subsubsection{
Translation to logics
}


\label{Section Reac-GenPref-Logic}
\index{Section Reac-GenPref-Logic}

We turn to the translation to logics.

\bp

$\hspace{0.01em}$


\label{Proposition Higher-Repr}

Let $ \xcn $ be a logic for $ \xdl.$ Set $T^{ \xdm }:=Th( \xbm_{ \xdm
}(M(T))),$ where $ \xdm $ is a
generalized preferential structure, and $ \xbm_{ \xdm }$ its choice
function. Then

(1) there is a level 2 preferential structure $ \xdm $ s.t. $ \ol{ \ol{T}
}=T^{ \xdm }$ iff
(LLE), (CCL), (SC) hold for all $T,T' \xcc \xdl.$

(2) there is a level 3 essentially smooth preferential structure
$ \xdm $ s.t. $ \ol{ \ol{T} }=T^{ \xdm }$ iff
(LLE), (CCL), (SC), $( \xcc \xcd )$ hold for all $T,T' \xcc \xdl.$

\ep

The proof is an immediate consequence of Corollary \ref{Corollary Eta-Rho}
(2),
Fact \ref{Fact X-Sub-X'},
Proposition \ref{Proposition Level-3-Repr},
and Proposition \ref{Proposition Alg-Log} (10) and (11).

We leave aside the generalization of preferential structures to attacking
structures relative to $ \xbh,$ as this can cause problems, without
giving real
insight: It might well be that $ \xbr (X) \xcC \xbh (X),$ but, still, $
\xbr (X)$ and $ \xbh (X)$ might
define the same theory - see Section \ref{Section Log-Introduction}.

\vspace{7mm}


\vspace{7mm}

\vspace{7mm}


\end{document}